\documentclass[12pt]{article}
\usepackage{amsmath,amssymb,amsthm,amscd,latexsym,}
 \font\goth=eusm10
\usepackage{graphicx}

\makeatletter
\renewcommand{\mod}[1]{\allowbreak \if@display \mkern 8mu \else
\mkern 5mu\fi {\operator@font mod}\,\,#1}
\makeatother

\newcommand{\iso}{\cong}
\newcommand{\ga}{\gamma}
\newcommand{\bc}{\mathbb C}
\newcommand{\bn}{\mathbb N}
 \newcommand{\bq}{\mathbb Q}
\newcommand{\br}{\mathbb R}
\newcommand{\bz}{\mathbb Z}
\newcommand{\bff}{\mathbb F}
\newcommand{\bl}{\mathbb L}
\newcommand{\bk}{\mathbb K}
\newcommand{\bo}{\mathbb O}
\newcommand{\bp}{\mathbb P}
\newcommand{\E}{\mathcal E}
\newcommand{\Oc}{\mathcal O}
\newcommand{\wf}{\widetilde{f}}
\newcommand{\wH}{\widetilde{H}}
\newcommand{\wth}{\widetilde{h}}
\newcommand{\wx}{\widetilde{x}}
\newcommand{\wy}{\widetilde{y}}
\newcommand{\aaa}{\mathbb A}
\newcommand{\ddd}{\mathbb D}
\newcommand{\eee}{\mathbb E}

\DeclareMathOperator{\Aut}{{\rm Aut}\,}
\DeclareMathOperator{\Hom}{Hom}
\DeclareMathOperator{\Mod}{Mod}
\DeclareMathOperator{\Pic}{Pic}
\DeclareMathOperator{\rk}{rk}
\DeclareMathOperator{\ch}{ch}
\DeclareMathOperator{\td}{td}
\DeclareMathOperator{\Tyu}{Tyu}
\DeclareMathOperator{\discr}{{\rm discr}\,}
\DeclareMathOperator{\cha}{{\rm char}\,}
\newtheorem{theorem}{Theorem}

\newcommand\F{\mathcal F}
\newcommand\Hh{\mathcal H}
\newcommand\N{\mathcal N}
\newcommand\La{\mathcal L}
\newcommand\Ka{{\mathcal K}}
\newcommand\G{\mathcal G}
\newcommand\M{\mathcal M}

\newcommand\SSS{\mathfrak S}
\newcommand\AAA{\mathfrak A}
\newcommand\DDD{\mathfrak D}

\begin{document}
\title{Classification of degenerations and Picard lattices of K\"ahlerian K3 surfaces with
the symplectic automorphism group $(C_2)^2$. }
\date{7 March 2021}
\author{Viacheslav V. Nikulin}
\maketitle

\begin{abstract} In our papers of 2013---2018, we classified
degenerations and Picard lattices of K\"ahle\-rian K3 surfaces with finite
symplectic automorphism groups of high order.
For remaining groups of small order: $D_6$, $C_4$, $(C_2)^2$, $C_3$, $C_2$ and $C_1$,
it was not completely considered because each of these cases requires very long and difficult considerations and
calculations. Cases of $D_6$ and $C_4$ were recently completely considered in \cite{Nik13} and \cite{Nik14}.

Here we consider the analogous complete classification for the group $(C_2)^2$ of the order $4$.
\end{abstract}

\centerline{Dedicated to the memory of Vitaliy Sergeevich Makarov}

\section{Introduction}
\label{sec:introduction}
In our papers \cite{Nik7} ---\cite{Nik12}, we considered classification of
degenerations and Picard lattices of K\"ahlerian K3 surfaces with finite
symplectic automorphism groups. In \cite{Nik7}---\cite{Nik12}, we completed this
classification for the groups of high order. For remaining groups of small orders:
$D_6$, $C_4$, $(C_2)^2$, $C_3$, $C_2$ and $C_1$, the complete classification
was not considered in these papers.
Here $D_n$ is the dihedral group, and $C_n$ is the cyclic group of order $n$.
For these groups, there are too many cases to consider, and it is better to consider each of
these groups separately. The cases of $D_6$ and $C_4$ were recently completely considered
in \cite{Nik13}, \cite{Nik14}.

In this paper, we consider similar complete classification for the group $(C_2^2)$ of the
order $4$.

This classification is given in Tables 1, 2 and 3 of Section \ref{sec2:degen}.

We use the same notations, methods and results as in our papers \cite{Nik7}---\cite{Nik13}.
See more details in Section \ref{sec2:degen}.

A preliminary variant of this paper was published in preprint \cite{Nik14}.


\section{Classification of degenerations and\\
Picard lattices $S$ of K\"ahlerian K3 surfaces with\\
finite symplectic automorphism groups.\\
General theory and the case of the group $(C_2)^2$.}
\label{sec2:degen}

Here we remind to a reader our methods and results from \cite{Nik7} --- \cite{Nik13}
for classification of degenerations and Picard lattices of K\"ahlerian K3 surfaces
with finite symplectic automorphism groups. Here we apply them to
the group $(C_2)^2$ of the order $4$
which was not completely considered in these papers.

Let $X$ be a K\"ahlerian K3 surface (e. g. see \cite{Sh}, \cite{PS},
\cite{BR}, \cite{Siu}, \cite{Tod}
about such surfaces). That is $X$ is a non-singular compact
complex surface with the trivial canonical class $K_X$, and
with the irregularity $q(X)=0$. It is K\"ahlerian by \cite{Siu} and \cite{Tod}.
Then $H^2(X,\bz)$ with the intersection pairing
is an even unimodular lattice $L_{K3}$
of the signature $(3,19)$. Further, a lattice means a non-degenerate integer
symmetric bilinear form of a finite rank.
For a non-zero holomorphic $2$-form $\omega_X\in \Omega^2[X]$,
we have $H^{2,0}(X)=\Omega^2[X]=\bc\omega_X$. The primitive sublattice
$$
S_X=H^2(X,\bz)\cap H^{1,1}(X)=\{x\in H^2(X,\bz)\ |\ x\cdot H^{2,0}(X)=0 \ \}\subset H^2(X,\bz)
$$
is the {\it Picard lattice} of $X$
generated by first Chern classes of all line bundles over $X$.
We remind that a primitive sublattice means that $H^2(X,\bz)/S_X$ has no torsion.

Let $G$ be a finite symplectic automorphism group of a K3 surface $X$.
See \cite{Nik-1/2}, \cite{Nik0}, \cite{Nik1}, \cite{Muk}, \cite{Kon}, \cite{Xiao}, \cite{Hash}
about their classification.
Here symplectic means
that for any $g\in G$, for a non-zero holomorphic $2$-form $\omega_X\in
H^{2,0}(X)=\Omega^2[X]=\bc\omega_X$, one has $g^\ast(\omega_X)=\omega_X$.
For an $G$-invariant sublattice $M\subset H^2(X,\bz)$, we denote by
$M^G=\{x\in M\ |\ G(x)=x\}$ the {\it fixed sublattice of $M$,} and by
$M_G=(M^G)^\perp_M$ the {\it coinvariant sublattice of $M$.}
By \cite{Nik-1/2}, \cite{Nik0}, the coinvariant lattice $S_G=H^2(X,\bz)_G=(S_X)_G$ is
{\it Leech type lattice:} it is negative definite,
it has no elements with square $(-2)$, $G$ acts trivially
on the discriminant group $A_{S_G}=(S_G)^\ast/S_G$, and $(S_G)^G=\{0\}$.
For a general pair $(X,G)$, the $S_G=S_X$, and non-general $(X,G)$
can be considered as K\"ahlerian K3 surfaces
with the condition $S_G\subset S_X$
on the Picard lattice (in terminology of \cite{Nik0}).
The dimension of their
moduli is equal to $20-\rk S_G$.

Let $E\subset X$ be a non-singular irreducible rational curve
(that is $E\cong \bp^1$).
It is equivalent to: $\alpha=cl(E)\in S_X$, $\alpha^2=-2$,
$\alpha$ is effective
and $\alpha$ is numerically effective: $\alpha\cdot D\ge 0$
for any irreducible curve
$D$ on $X$ such that $cl(D)\not=\alpha$.

Let us consider $t$ non-singular irreducible rational curves $E_1,\dots, E_t$ on $X$
with classes $\alpha_i=cl(E_i)\in S_X$ such that
their orbits $G(E_1),\dots, G(E_t)$ are different.
Let us consider the primitive sublattice
$$
S=[S_G, G(\alpha_1),\dots, G(\alpha_t)]_{pr}=
[S_G, \alpha_1,\dots, \alpha_t]_{pr}\subset S_X
$$
of $S_X$ generated by the coinvariant sublattice $S_G$ and
all classes of the orbits $G(E_1),\dots G(E_t)$ and {\it assume that $S$ is negative
definite.}  Since $S_G$ has no elements with square $(-2)$, it follows
that $\rk S=\rk S_G+t$ and elements of orbits $G(\alpha_1),\dots G(\alpha_t)$
define the basis of the root system $\Delta(S)$ of all elements with square $(-2)$ of $S$.
All curves $G(E_1),\dots, G(E_t)$ of $X$
can be contracted to Du Val singularities of types of
connected components of the Dynkin diagram of the basis. The group $G$
will act on the corresponding singular K3 surface
$\overline{X}$ with these Du Val
singularities. For a general such collection $(X,G,G(E_1),\dots G(E_t))$,
the Picard lattice $S_X=S$,
and such collection can be considered as {\it a degeneration
of codimension $t$} of
K\"ahlerian K3 surfaces $(X,G)$ with the finite symplectic
automorphism group $G$.
Really, the dimension of moduli of K\"ahlerian K3 surfaces
with the condition $S\subset S_X$
on the Picard lattice is equal to $20-\rk S=20-\rk S_G-t$.

We can consider only maximal finite symplectic
automorphism groups $G$ with the
same coinvariant lattices $S_G$, that is $G=Clos(G)$.
By Global Torelli Theorem for K3 surfaces \cite{PS}, \cite{BR},  this
is equivalent to
$$
G|S_G=\{ g\in O(S_G)\ |\ g\ is\ identity\ on\ A_{S_G}=(S_G)^\ast/S_G \}.
$$

The {\it type of the degeneration} is given by the Dynkin
diagrams and subdiagrams
$$
(Dyn(G(\alpha_1)),\dots, Dyn(G(\alpha_t)))\subset
Dyn(G(\alpha_1)\cup\dots \cup G(\alpha_t))
$$
and their types. Numberings of subdiagrams $Dyn(G(\alpha_i))$
and connected components  of $Dyn(G(\alpha_1)\cup\dots \cup G(\alpha_t))$ must agree.
In difficult cases, we also consider the matrix of subdiagrams which is defined by
$$
(Dyn (G(\alpha_i)),\, Dyn ( G(\alpha_j)))\subset Dyn(G(\alpha_i)\cup G(\alpha_j))
$$
and their types for $1\le i<j\le t$.

By Global Torelli Theorem for K3 surfaces \cite{PS},
\cite{BR}, the type of the abstract group $G=Clos(G)$
which is equivalent
to the isomorphism class of the coinvariant lattice $S_G$,
and the type of the degeneration give the main invariants of
the degeneration.

For groups $G=Clos(G)$ of order $|G|>4$ and $G=C_4$ classification of degenerations of arbitrary
codimension and Picard lattices $S$ is given in tables of \cite{Nik9}---\cite{Nik13}.
In this paper, applying the same results and methods, we give similar classification of
degenerations of arbitrary codimension and Picard lattices $S$ for
the group $G=(C_2)^2$. It is given in Table 1 where we give all possible types of
degenerations in the column $Deg$. We use the standard notations $\aaa_n$, $\ddd_n$, $\eee_n$
for connected components of Dynkin diagrams with $n$ vertices (see \cite{Bou}), and $A_n$, $D_n$, $E_n$
the corresponding root lattices. By $\Gamma_1\amalg\Gamma_2$ we denote the disjoint union
of diagrams $\Gamma_1$ and $\Gamma_2$ and by $\oplus$ the orthogonal sum of lattices.
By $kM$ we denote the orthogonal sum of $k$ copies of a lattice $M$.
For a given degeneration $Deg$,
at the same line, using results of \cite{Nik1}, we give the genus of the lattices $S$ by
$\rk S$ and the discriminant quadratic form $q_S$ in notations
of \cite{CS}. By $\ast$, we denote cases
when we also prove that
the lattice $S$ is unique up to isomorphisms for the given type.

In Table 2, we give some markings by Niemeier lattices $N_j$, $j=1,\dots, 24$,
of lattices $S$ of a degeneration $Deg$ together with the action of $G$ on $S$
and its orbits $G(\alpha_1),\dots G(\alpha_t)$
in notations of \cite{Nik8}, \cite{Nik9}. For these markings, $G$ is
denoted by $H$. They give a description of possible $(S,G,G(\alpha_1),\dots G(\alpha_t))$
by integer lattices. For these markings, we denote by $\ast$ markings
which imply uniqueness of the lattice $S$, up to isomorphisms,
where we use ideas by Hashimoto \cite{Hash} and our methods in \cite{Nik9}---\cite{Nik12}.
Differently from \cite{Nik9}---\cite{Nik12}, we give only part of analogous to
\cite{Nik9}---\cite{Nik12} markings by Niemeier lattices, which is sufficient for
the description  of lattices $S$ of degenerations $Deg$ (otherwise, there will be too
many markings and cases). Moreover (analogous to methods which we discussed in
\cite{Nik9}---\cite{Nik12}), we also use degenerations of smaller codimensions and
their markings by Niemeier lattices if the considered degeneration can be reduced to
them because of calculations of geneses of lattices $S$ of the Table 1. This significantly
reduces the number of considered cases.

Additionally to Tables 1 and 2, in Table 3 we give the List 1 which is important for
classification of Picard lattices $S$ of K3 surfaces (see \cite{Nik12}).
In Table 1, we denote by $o$ (old) cases when the
degeneration of K3 surfaces with symplectic automorphism
group $(C_2)^2$ has, actually,  the maximal finite symplectic automorphism group $G$ which
contains $(C_2)^2$ and $|G|>|(C_2)^2|=4$. The group $G$ has less orbits and less codimension
of the degeneration than $(C_2)^2$. For classification of Picard lattices $S$
of K3 surfaces with maximal finite symplectic automorphism group $(C_2)^2$, lines of Table 1
which are denoted by $o$ must be removed.
In Table 3, we give the List 1 of all these cases (which is similar
to the List 1 of \cite{Nik12}) where the degeneration of the group $(C_2)^2$ is shown to the left,
and the group $G$ (defined by ${\bf n}$ as an abstract group) and its degeneration
(classified in \cite{Nik10}---\cite{Nik13}) is shown to the right from the sign $\Longleftarrow$.
Cases of Table 1 which are marked by $o$ are also very interesting, they are similar to the case of Kummer
surfaces which give the degeneration of type $16\aaa_1$ of K3 surfaces
with zero Picard lattice and trivial symplectic automorphism group, but the symplectic automorphism
group of a general Kummer surface is $(C_2)^4$. See \cite{Nik-1}, \cite{Nik12}.



\begin{table}
\label{table1}
\caption{Some types and lattices $S$ of degenerations of K\"ahlerian K3 surfaces
with the symplectic automorphism group $G=(C_2)^2$.}



\begin{tabular}{|c||c|c|c|c|c|c|c|c|}
\hline
 {\bf n}& $|G|$& $i$&  $G$  & $\rk S_G$ &$q_{S_G}$  &$Deg$& $\rk S$ &$q_S$ \\
\hline
\hline
$3$   &$4$& $2$ &$(C_2)^2$  &   $12$     &$2_{II}^{-6},4_{II}^{-2}$ & & $12$ & $2_{II}^{-6},4_{II}^{-2}$ $\ast$\\
\hline
\hline
 $3$  &$4$& $2$ & $(C_2)^2$ & $12$      &$2_{II}^{-6},4_{II}^{-2}$    &  $\aaa_1$   & $13$  &
$2_3^{+7},4_{II}^{+2}$ $\ast$ \\
\hline
      &   &     &       &           &                        &  $2\aaa_1$  & $13$  &
$2_{II}^{-4},4_{7}^{-3}$ $\ast$ \\
\hline
      &   &     &       &           &                        &  $4\aaa_1$  & $13$  &
$2_{II}^{-6},8_{3}^{-1}$ $\ast$ \\
\hline
\hline
 & &  & &  & &  $(\aaa_1,\aaa_1)\subset 2\aaa_1$  & $14$ &$2_{6}^{-6},4_{II}^{+2}$   \\
\hline
 & &  & &  & &  $(\aaa_1,2\aaa_1)\subset 3\aaa_1$  & $14$ &$2_{1}^{+5},4_{1}^{+3}$\\
\hline
 & &  & &  & &  $(\aaa_1,2\aaa_1)\subset \aaa_3$  & $14$ &$2_{6}^{-6},4_{II}^{+2}$\\
\hline

\end{tabular}
\end{table}

\begin{table}
\begin{tabular}{|c||c|c|c|c|c|c|c|c|}
\hline
 {\bf n}& $G$ & $\rk S_G$ & $Deg$ & $\rk S$ & $q_S$ \\
\hline
\hline
 $3$  & $(C_2)^2$ & $12$ &
 $(\aaa_1,2\aaa_1,4\aaa_1)\subset \aaa_3\amalg 4\aaa_1$  &
 $15$ &$2_6^{-6},8_7^{+1}$ \\
\hline
&  & &  $(\aaa_1,4\aaa_1,4\aaa_1)\subset 9\aaa_1$  & $15$ &$2_{3}^{-5},4_2^{+2}$ \\
\hline
&  & &  $(\aaa_1,4\aaa_1,4\aaa_1)\subset \aaa_1\amalg 4\aaa_2$
 & $15$ &$2_{3}^{+7},3^{+1}$ \\
\hline
& & &
 $\left(\begin{array}{rrr}
2\aaa_1 & (4\aaa_1)_{I} & (4\aaa_1)_I  \\
        & 2\aaa_1     & (4\aaa_1)_I   \\
        &             & 2\aaa_1
\end{array}\right)\subset 6\aaa_1$
 & $15$ &$4_{1}^{+5}$ \\
\hline
 & &  &
$\left(\begin{array}{rrr}
2\aaa_1 & (4\aaa_1)_{II} & (4\aaa_1)_I  \\
        & 2\aaa_1     & (4\aaa_1)_I   \\
        &             & 2\aaa_1
\end{array}\right)\subset 6\aaa_1$
 & $15$ &$2_{II}^{-4},4_{5}^{-3}$ \\
\hline
& &  &   $(2\aaa_1,2\aaa_1,2\aaa_1)\subset 2\aaa_1\amalg 2\aaa_2$ &
$15$ &$2_{II}^{-2},4_{3}^{+3},3^{-1}$ \\
\hline
 & &  &
 $\left(\begin{array}{rrr}
 2\aaa_1 & 2\aaa_2 & (4\aaa_1)_{II}  \\
         & 2\aaa_1 &  2\aaa_2   \\
         &         & 2\aaa_1
\end{array}\right)\subset 2\aaa_3$
& $15$ &$2_{II}^{+6},8_1^{+1}$ \\
\hline

\end{tabular}
\end{table}

\begin{table}
\begin{tabular}{|c||c|c|c|c|c|c|c|c|}
\hline
 {\bf n}& $G$ & $\rk S_G$ & $Deg$& $\rk S$ & $q_S$ \\
\hline
\hline
 $3$  & $(C_2)^2$ & $12$ &
 $(\aaa_1,\aaa_1,\aaa_1,4\aaa_1)\subset 7\aaa_1$  & $16$ &$2_{1}^{+5},8_7^{+1}$ \\
\hline
& &  &
 $(\aaa_1,\aaa_1,(2\aaa_1,2\aaa_1)_I)\subset 6\aaa_1$  & $16$ &$2_{4}^{-2},4_0^{-4}$ \\
 \hline
 & & &
   $\left(\begin{array}{rrrr}
 \aaa_1 & 2\aaa_1 & 3\aaa_1  & 3\aaa_1  \\
         & \aaa_1 & \aaa_3   & 3\aaa_1   \\
         &         & 2\aaa_1 & (4\aaa_1)_I \\
         &         &         & 2\aaa_1
\end{array}\right)\subset \aaa_1\amalg \aaa_3\amalg 2\aaa_1$
 & $16$ & $2_3^{-3},4_{1}^{+3}$ \\
 \hline
&  &  &
   $\left(\begin{array}{rrrr}
 \aaa_1  & \aaa_3  & 2\aaa_1  & 3\aaa_1  \\
         & 2\aaa_1 & 3\aaa_1  & (4\aaa_1)_I   \\
         &         & \aaa_1   & \aaa_3 \\
         &         &          & 2\aaa_1
\end{array}\right)$
 & $16$ & $2_0^{+4},4_{II}^{+2}$ \\
& & & $\subset \aaa_3\amalg \aaa_3$ & &\\
 \hline
&   &  &
   $\left(\begin{array}{rrrr}
 \aaa_1  & \aaa_3  & 2\aaa_1  & 3\aaa_1  \\
         & 2\aaa_1 & 3\aaa_1  & (4\aaa_1)_{II}   \\
         &         & \aaa_1   & \aaa_3 \\
         &         &          & 2\aaa_1
\end{array}\right)$
 & $16$ & $2_{II}^{+4},4_{II}^{+2}$ \\
& & & $\subset \aaa_3\amalg \aaa_3$ & & \\
\hline

\end{tabular}
\end{table}

\begin{table}
\begin{tabular}{|c||c|c|c|c|c|c|c|c|}
\hline
 {\bf n}& $G$ & $\rk S_G$ & $Deg$ & $\rk S$ & $q_S$ \\
\hline
 $3$  & $(C_2)^2$ & $12$ &
 $\left(\begin{array}{rrrr}
 \aaa_1 & \aaa_3 & 3\aaa_1  & 3\aaa_1  \\
         & 2\aaa_1 & (4\aaa_1)_{II} & (4\aaa_1)_I   \\
         &         & 2\aaa_1 & (4\aaa_1)_I   \\
         &         &         & 2\aaa_1
\end{array}\right)\subset \aaa_3\amalg 4\aaa_1$
 & $16$ &$2_{II}^{+2},4_{0}^{+4}$ \\
 \hline
     & & &
 $\left(\begin{array}{rrrr}
 \aaa_1 & \aaa_3 & 3\aaa_1  & 3\aaa_1  \\
         & 2\aaa_1 & (4\aaa_1)_{I} & (4\aaa_1)_I   \\
         &         & 2\aaa_1 & 2\aaa_2   \\
         &         &         & 2\aaa_1
\end{array}\right)\subset \aaa_3\amalg 2\aaa_2$
 & $16$ &$2_{6}^{+4},4_{II}^{+2},3^{-1}$ \\
\hline
& & &
 $\left(\begin{array}{rrrr}
 2\aaa_1 & (4\aaa_1)_{I}  & (4\aaa_1)_{I}  & 3\aaa_1  \\
         & 2\aaa_1 & 2\aaa_2 & 3\aaa_1   \\
         &         & 2\aaa_1 & \aaa_3   \\
         &         &         & \aaa_1
\end{array}\right)\subset 2\aaa_1\amalg \aaa_5$
 & $16$ &$2_{5}^{+3},4_{3}^{-3}$ \\
\hline
& & &
 $\left(\begin{array}{rrrr}
 2\aaa_1 & 2\aaa_2 & (4\aaa_1)_{II}  & 3\aaa_1  \\
         & 2\aaa_1 & 2\aaa_2         & 3\aaa_1   \\
         &         & 2\aaa_1         & \aaa_3   \\
         &         &                 & \aaa_1
\end{array}\right)\subset \aaa_7$
 & $16$ &$2_{II}^{+4},4_{II}^{+2}$ \\
\hline
 & & &
 $(\aaa_1,(2\aaa_1,2\aaa_1)_I,4\aaa_1)\subset 9\aaa_1$  & $16$ &$2_7^{+3},4_2^{+2},8_{3}^{-1}$ \\
\hline
 & & &
 $(\aaa_1,2\aaa_1,2\aaa_1,4\aaa_1)\subset \aaa_1\amalg 2\aaa_2\amalg 4\aaa_1$ & $16$ &$2_3^{-5},8_{3}^{-1},3^{-1}$ \\
\hline
& & &
 $\left(\begin{array}{rrrr}
 \aaa_1 & 3\aaa_1 & 3\aaa_1  & 5\aaa_1  \\
         & 2\aaa_1 & (4\aaa_1)_{I} & 6\aaa_1   \\
         &         & 2\aaa_1 & 2\aaa_3   \\
         &         &         & 4\aaa_1
\end{array}\right)\subset 3\aaa_1\amalg 2\aaa_3$
 & $16$ &$2_{3}^{-3},4_{1}^{+3}$ \\
\hline

\end{tabular}
\end{table}

\begin{table}
\begin{tabular}{|c||c|c|c|c|c|c|c|c|}
\hline
 {\bf n}& $G$ &$Deg$& $\rk S$ &$q_S$ \\
\hline
\hline
 $3$  & $(C_2)^2$ &
     $\left(\begin{array}{rrrr}
 \aaa_1 & \aaa_3 & 3\aaa_1  & 5\aaa_1  \\
         & 2\aaa_1 & (4\aaa_1)_{I} & 6\aaa_1   \\
         &         & 2\aaa_1 & 2\aaa_3   \\
         &         &         & 4\aaa_1
\end{array}\right)\subset \aaa_3\amalg 2\aaa_3$
 & $16$ &$2_{4}^{-4},4_{0}^{+2}$   \\
\hline
& &
     $\left(\begin{array}{rrrr}
 \aaa_1 & \aaa_3 & 3\aaa_1  & 5\aaa_1  \\
         & 2\aaa_1 & (4\aaa_1)_{II} & 6\aaa_1   \\
         &         & 2\aaa_1 & 2\aaa_3   \\
         &         &         & 4\aaa_1
\end{array}\right)\subset \aaa_3\amalg 2\aaa_3$
 & $16$ &$2_{II}^{+4},4_{0}^{+2}$   \\
\hline
& &
     $\left(\begin{array}{rrrr}
 4\aaa_1 & 2\aaa_3 & 6\aaa_1  & 5\aaa_1  \\
         & 2\aaa_1 & 2\aaa_2 & 3\aaa_1   \\
         &         & 2\aaa_1 & 3\aaa_1   \\
         &         &         & \aaa_1
\end{array}\right)\subset 2\ddd_4\amalg \aaa_1$
 & $16$ &$2_{7}^{+5},8_{1}^{+1}$   \\
\hline
 & &
 $(4\aaa_1,2\aaa_1,2\aaa_1,\aaa_1)\subset 4\aaa_1\amalg \aaa_5$ & $16$ &$2_5^{-5},8_{7}^{+1}$ \\
\hline
 & &
 $(\aaa_1,2\aaa_1,4\aaa_1,4\aaa_1)\subset 11\aaa_1$ & $16$ &$2_3^{-3},4_{1}^{+3}$ \\
\hline
 & &
 $(\aaa_1,2\aaa_1,4\aaa_1,4\aaa_1)\subset 3\aaa_1\amalg 4\aaa_2$ & $16$ &$2_5^{+5},4_{5}^{-1},3^{+1}$ \\
\hline
 & &
 $(\aaa_1,2\aaa_1,4\aaa_1,4\aaa_1)\subset \aaa_3\amalg 4\aaa_2$ & $16$ &$2_2^{+6},3^{+1}$ \\
\hline

\end{tabular}
\end{table}

\begin{table}
\begin{tabular}{|c||c|c|c|c|c|c|c|c|}
\hline
 {\bf n}& $G$  & $Deg$& $\rk S$ &$q_S$ \\
\hline
\hline
 $3$  & $(C_2)^2$ &
     $\left(\begin{array}{rrrr}
 2\aaa_1 & (4\aaa_1)_{I} & (4\aaa_1)_I  & (4\aaa_1)_I  \\
         & 2\aaa_1 & 2\aaa_2 & (4\aaa_1)_I   \\
         &         & 2\aaa_1 & (4\aaa_1)_{I}   \\
         &         &         & 2\aaa_1
\end{array}\right)\subset 2\aaa_1\amalg 2\aaa_2\amalg 2\aaa_1$
 & $16$ &$4_{6}^{-4},3^{-1}$   \\
\hline
 & &
       $\left(\begin{array}{rrrr}
 2\aaa_1 & (4\aaa_1)_{I} & (4\aaa_1)_I  & (4\aaa_1)_{II}  \\
         & 2\aaa_1 & 2\aaa_2 & (4\aaa_1)_I   \\
         &         & 2\aaa_1 & (4\aaa_1)_{I}   \\
         &         &         & 2\aaa_1
\end{array}\right)\subset 2\aaa_1\amalg 2\aaa_2\amalg 2\aaa_1$
 & $16$ &$2_{II}^{-4},4_2^{+2},3^{-1}$   \\
\hline
 & &
 $(2\aaa_1,2\aaa_1,2\aaa_1,2\aaa_1)\subset 2\aaa_2\amalg 2\aaa_2$ & $16$ &$2_{II}^{-2},4_{II}^{-2},3^{+2}$ \\
\hline
& &
     $\left(\begin{array}{rrrr}
 2\aaa_1 & (4\aaa_1)_{I} & (4\aaa_1)_I  & (4\aaa_1)_I  \\
         & 2\aaa_1 & 2\aaa_2 & (4\aaa_1)_{II}   \\
         &         & 2\aaa_1 & 2\aaa_2   \\
         &         &         & 2\aaa_1
\end{array}\right)\subset 2\aaa_1\amalg 2\aaa_3$
 & $16$ &$2_{II}^{+4},4_{7}^{+1},8_1^{+1}$   \\
\hline
& &
   $\left(\begin{array}{rrrr}
 2\aaa_1 & (4\aaa_1)_{I} & (4\aaa_1)_I  & 6\aaa_1  \\
         & 2\aaa_1 & (4\aaa_1)_I & 6\aaa_1   \\
         &         & 2\aaa_1 & 6\aaa_1   \\
         &         &         & 4\aaa_1
\end{array}\right)\subset 10\aaa_1$
 & $16$ &$4_{7}^{-3},8_5^{-1}$   \\
\hline
& &
   $\left(\begin{array}{rrrr}
 2\aaa_1 & (4\aaa_1)_{II} & (4\aaa_1)_I  & 6\aaa_1  \\
         & 2\aaa_1 & (4\aaa_1)_I & 6\aaa_1   \\
         &         & 2\aaa_1 & 6\aaa_1   \\
         &         &         & 4\aaa_1
\end{array}\right)\subset 10\aaa_1$
 & $16$ &$2_{II}^{-4},4_{3}^{-1},8_1^{+1}$   \\
\hline

\end{tabular}
\end{table}

\begin{table}
\begin{tabular}{|c||c|c|c|c|c|c|c|c|}
\hline
 {\bf n}& $G$  & $Deg$& $\rk S$ &$q_S$ \\
\hline
 $3$  & $(C_2)^2$ &
      $\left(\begin{array}{rrrr}
 4\aaa_1 & 2\aaa_3 & 6\aaa_1  & 6\aaa_1  \\
         & 2\aaa_1 & (4\aaa_1)_{I} & (4\aaa_1)_I   \\
         &         & 2\aaa_1 & 2\aaa_2   \\
         &         &         & 2\aaa_1
\end{array}\right)\subset 2\aaa_3\amalg 2\aaa_2$
 & $16$ &  $2_{II}^{-2},4_{2}^{+2},3^{-1}$   \\
\hline
& &
      $\left(\begin{array}{rrrr}
 4\aaa_1 & 6\aaa_1 & 6\aaa_1  & 6\aaa_1  \\
         & 2\aaa_1 & 2\aaa_2 & (4\aaa_1)_I   \\
         &         & 2\aaa_1 & (4\aaa_1)_I   \\
         &         &         & 2\aaa_1
\end{array}\right)\subset 4\aaa_1\amalg 2\aaa_2\amalg 2\aaa_1$
 & $16$ &  $2_{II}^{-2},4_{3}^{-1},8_3^{-1},3^{-1}$   \\
\hline
& &
      $\left(\begin{array}{rrrr}
 4\aaa_1 & 6\aaa_1 & 6\aaa_1  & 6\aaa_1  \\
         & 2\aaa_1 & 2\aaa_2 & (4\aaa_1)_{II}   \\
         &         & 2\aaa_1 & 2\aaa_2   \\
         &         &         & 2\aaa_1
\end{array}\right)\subset 4\aaa_1\amalg 2\aaa_3$
 & $16$ &  $2_{II}^{+4},4_{II}^{+2}$   \\
\hline
& &
      $\left(\begin{array}{rrrr}
 4\aaa_1 & 2\aaa_3 & 6\aaa_1  & 6\aaa_1  \\
         & 2\aaa_1 & 2\aaa_2 & (4\aaa_1)_{I}   \\
         &         & 2\aaa_1 & (4\aaa_1)_I   \\
         &         &         & 2\aaa_1
\end{array}\right)\subset 2\ddd_4\amalg 2\aaa_1$
 & $16$ &  $2_{II}^{-2},4_{3}^{-1},8_1^{+1}$   \\
\hline
& &
$(4\aaa_1,2\aaa_1,2\aaa_1,2\aaa_1)\subset 2\ddd_5$  &
$16$ &$2_{II}^{+6}$ \\
\hline

& &
$((2\aaa_1,2\aaa_1)_I,4\aaa_1,4\aaa_1)\subset 12\aaa_1$  &
$16$ &$4_{0}^{+4}$ \\
\hline
& &
$((2\aaa_1,2\aaa_1)_{II},4\aaa_1,4\aaa_1)\subset 12\aaa_1$  &
$16$ &$2_{II}^{+4},4_{II}^{+2}$ \\
\hline

\end{tabular}
\end{table}

\begin{table}
\begin{tabular}{|c||c|c|c|c|c|c|c|c|}
\hline
 {\bf n}& $G$  & $Deg$& $\rk S$ &$q_S$ \\
\hline
\hline
 $3$  & $(C_2)^2$ &
       $\left(\begin{array}{rrrr}
 4\aaa_1 & 2\aaa_3 & 6\aaa_1  & 8\aaa_1  \\
         & 2\aaa_1 & 2\aaa_2 & 6\aaa_1   \\
         &         & 2\aaa_1 & 6\aaa_1   \\
         &         &         & 4\aaa_1
\end{array}\right)\subset 2\ddd_4\amalg 4\aaa_1$
 & $16$ &  $2_{II}^{+2},4_{II}^{+2}$   \\
\hline
&  &
       $\left(\begin{array}{rrrr}
 4\aaa_1 & 4\aaa_2 & 6\aaa_1  & 6\aaa_1  \\
         & 4\aaa_1 & 2\aaa_3 & 6\aaa_1   \\
         &         & 2\aaa_1 & (4\aaa_1)_I   \\
         &         &         & 2\aaa_1
\end{array}\right)\subset 2\aaa_5\amalg 2\aaa_1$
 & $16$ &  $2_{II}^{+2},4_{0}^{+2}$   \\
\hline
&  &
       $\left(\begin{array}{rrrr}
 4\aaa_1 & 4\aaa_2 & 6\aaa_1  & 6\aaa_1  \\
         & 4\aaa_1 & 2\aaa_3 & 6\aaa_1   \\
         &         & 2\aaa_1 & (4\aaa_1)_{II}   \\
         &         &         & 2\aaa_1
\end{array}\right)\subset 2\aaa_5\amalg 2\aaa_1$
 & $16$ &  $2_{II}^{+6}$   \\
\hline
& &
$(4\aaa_1,4\aaa_1,2\aaa_1,2\aaa_1)\subset 2\eee_6$  &
$16$ &$2_{II}^{+4}$ \\
\hline
& &
$(2\aaa_1,4\aaa_1,4\aaa_1,4\aaa_1)\subset 2\aaa_3\amalg 8\aaa_1$  &
$16$ &$2_{II}^{+2},4_{II}^{+2}$ \\
\hline
& &
$(2\aaa_1,4\aaa_1,4\aaa_1,4\aaa_1)\subset 2\aaa_1\amalg 4\aaa_3$  &
$16$ &$2_{II}^{-2},4_{4}^{-2}$ \\
\hline
&  &
       $\left(\begin{array}{rrrr}
 4\aaa_1 & 4\aaa_2 & 8\aaa_1  & 6\aaa_1  \\
         & 4\aaa_1 & 4\aaa_2 & 6\aaa_1   \\
         &         & 4\aaa_1 & 2\aaa_3   \\
         &         &         & 2\aaa_1
\end{array}\right)\subset 2\aaa_7$
 & $16$ &  $2_{II}^{+4}$   \\
\hline
\hline

\end{tabular}
\end{table}

\begin{table}
\begin{tabular}{|c||c|c|c|c|c|c|c|c|}
\hline
 {\bf n}& $G$  & $Deg$& $\rk S$ &$q_S$ \\
\hline
\hline
 $3$  & $(C_2)^2$ &
$(\aaa_1,\aaa_1,\aaa_1,\aaa_1,4\aaa_1)\subset 8\aaa_1$  &
$17$ &$2_{II}^{+4},8_7^{+1}$  \\
\hline
&  &
       $\left(\begin{array}{rrrrr}
2\aaa_1 & \aaa_3 & (4\aaa_1)_{I}  &  3\aaa_1   &  3\aaa_1  \\
       & \aaa_1  & 3\aaa_1  &  2\aaa_1  &  2\aaa_1  \\
       &         & 2\aaa_1  & \aaa_3   &  3\aaa_1  \\
       &         &          & \aaa_1    & 2\aaa_1 \\
       &         &          &            & \aaa_1
\end{array}\right)$
 & $17$ & $2_{3}^{-3},4_{II}^{-2}$  \\
 &      & $\subset \aaa_3\amalg \aaa_3\amalg \aaa_1 $ & &    \\
 \hline
 & &
$(\aaa_1,\aaa_1,\aaa_1,2\aaa_1,4\aaa_1)\subset 2\aaa_1\amalg \aaa_3\amalg 4\aaa_1$  &
$17$ &$2_{4}^{-4},8_3^{-1}$  \\
\hline
 & &
$(\aaa_1,\aaa_1,\aaa_1,4\aaa_1,4\aaa_1)\subset 11\aaa_1$  &
$17$ &$2_{5}^{-3},4_6^{+2}$  \\
\hline
 & &
$(\aaa_1,\aaa_1,\aaa_1,4\aaa_1,4\aaa_1)\subset 3\aaa_1\amalg 4\aaa_2$  &
$17$ & $2_{1}^{+5},3^{+1}$  \\
\hline
 & &
     $\left(\begin{array}{rrrrr}
\aaa_1 & 2\aaa_1 & 3\aaa_1  &  3\aaa_1   &  3\aaa_1  \\
       & \aaa_1  & 3\aaa_1  & 3\aaa_1   &  3\aaa_1  \\
       &         & 2\aaa_1  & (4\aaa_1)_I&  (4\aaa_1)_I  \\
       &         &          & 2\aaa_1    &  (4\aaa_1)_I  \\
       &         &          &            &   2\aaa_1
\end{array}\right)\subset 8\aaa_1$
 & $17$ &$4_{7}^{+5}$   \\
\hline

\end{tabular}
\end{table}

\begin{table}
\begin{tabular}{|c||c|c|c|c|c|c|c|c|}
\hline
 {\bf n}& $G$  & $Deg$& $\rk S$ &$q_S$ \\
\hline
\hline
 $3$  & $(C_2)^2$ &
$\left(\begin{array}{rrrrr}
2\aaa_1 & (4\aaa_1)_I  & 3\aaa_1  &  (4\aaa_1)_I        &  3\aaa_1  \\
        & 2\aaa_1      & \aaa_3   &   (4\aaa_1)_{II}    &  3\aaa_1  \\
        &              & \aaa_1   &   3\aaa_1           &  2\aaa_1  \\
        &              &          & 2\aaa_1             &  \aaa_3 \\
       &               &          &                     &  \aaa_1
\end{array}\right)$
 & $17$ & $2_{II}^{+2},4_{7}^{+3}$  \\
 &      & $\subset 2\aaa_1\amalg \aaa_3\amalg \aaa_3$ & &    \\
 \hline
&  &
$\left(\begin{array}{rrrrr}
2\aaa_1 & (4\aaa_1)_I  & (4\aaa_1)_I        &  3\aaa_1   &  3\aaa_1 \\
        & 2\aaa_1      & 2\aaa_2            &  3\aaa_1   &  3\aaa_1  \\
        &              & 2\aaa_1            &  \aaa_3    &  3\aaa_1  \\
        &              &                    &  \aaa_1    &  2\aaa_1 \\
        &              &                    &            &  \aaa_1
\end{array}\right)$
 & $17$ & $2_{0}^{+2},4_{7}^{-3}$  \\
 &      & $\subset 2\aaa_1\amalg \aaa_5\amalg \aaa_1$ & &    \\
 \hline
  &  &
$\left(\begin{array}{rrrrr}
2\aaa_1 & 2\aaa_2      & 3\aaa_1            &  (4\aaa_1)_I   &  3\aaa_1 \\
        & 2\aaa_1      & \aaa_3             &  (4\aaa_1)_I   &  3\aaa_1  \\
        &              & \aaa_1             &  3\aaa_3       &  3\aaa_1  \\
        &              &                    &  2\aaa_1       &  \aaa_3 \\
        &              &                    &                &  \aaa_1
\end{array}\right)$
 & $17$ & $2_{3}^{-3},4_{II}^{+2}$  \\
 &      & $\subset \aaa_5\amalg \aaa_3$  & &    \\
 \hline

&  &
$(\aaa_1,\aaa_1,(2\aaa_1,2\aaa_1)_{I},4\aaa_1)\subset 10\aaa_1$  &
$17$ &$2_{6}^{-2},4_{6}^{-2},8_3^{-1}$  \\
\hline

\end{tabular}
\end{table}

\begin{table}
\begin{tabular}{|c||c|c|c|c|c|c|c|c|}
\hline
 {\bf n}& $G$  & $Deg$& $\rk S$ &$q_S$ \\
\hline
\hline
 $3$  & $(C_2)^2$ &
     $\left(\begin{array}{rrrrr}
\aaa_1 & \aaa_3   & 2\aaa_1  &  3\aaa_1          &  5\aaa_1  \\
       & 2\aaa_1  & 3\aaa_1  &  (4\aaa_1)_{II}   &  6\aaa_1  \\
       &         & \aaa_1    & \aaa_3            &  5\aaa_1  \\
       &         &           & 2\aaa_1           &  3\aaa_1  \\
       &         &           &                   &  4\aaa_1
\end{array}\right)$
 & $17$ &$2_{II}^{+4},8_{7}^{+1}$   \\
 & & $\subset \aaa_3\amalg \aaa_3\amalg 4\aaa_1$ & & \\
 \hline
& &
     $\left(\begin{array}{rrrrr}
\aaa_1 & 2\aaa_1 & 3\aaa_1  &  3\aaa_1       &  5\aaa_1  \\
       & \aaa_1  & 3\aaa_1  &  3\aaa_1       &  5\aaa_1  \\
       &         & 2\aaa_1  &  (4\aaa_1)_I   &  6\aaa_1  \\
       &         &          &  2\aaa_1    &  2\aaa_3 \\
       &         &          &             &  4\aaa_1
\end{array}\right)\subset 4\aaa_1\amalg 2\aaa_3$
 & $17$ & $2_6^{+2},4_{1}^{+3}$   \\
\hline
& &
     $\left(\begin{array}{rrrrr}
\aaa_1 & 2\aaa_1 & 3\aaa_1  &  3\aaa_1           &      5\aaa_1  \\
       & \aaa_1  & \aaa_3   &  3\aaa_1           &      5\aaa_1  \\
       &         & 2\aaa_1  &  (4\aaa_1)_I       &      6\aaa_1  \\
       &         &          & 2\aaa_1            &      2\aaa_3  \\
       &         &          &                    &      4\aaa_1
\end{array}\right) $
 & $17$ & $2_7^{+3},4_{0}^{+2}$     \\
 & & $\subset \aaa_1\amalg \aaa_3\amalg 2\aaa_3$ & & \\
\hline
 & &
$(\aaa_1,2\aaa_1,2\aaa_1,\aaa_1,4\aaa_1)\subset \aaa_1\amalg \aaa_5\amalg 4\aaa_1$  &
$17$ &$2_{0}^{+4},8_{7}^{+1}$     \\
\hline
 & &
$(\aaa_1,\aaa_1,2\aaa_1,4\aaa_1,4\aaa_1)\subset 12\aaa_1$  &
$17$ &$2_{6}^{+2},4_{1}^{+3}$     \\
\hline

\end{tabular}
\end{table}

\begin{table}
\begin{tabular}{|c||c|c|c|c|c|c|c|c|}
\hline
 {\bf n}& $G$  & $Deg$& $\rk S$ &$q_S$ \\
\hline
\hline
 $3$  & $(C_2)^2$ &
     $\left(\begin{array}{rrrrr}
2\aaa_1 & 2\aaa_2  & (4\aaa_1)_I  &  (4\aaa_1)_I   &  3\aaa_1  \\
        & 2\aaa_1  & (4\aaa_1)_I  &  (4\aaa_1)_I   &  3\aaa_1  \\
        &          & 2\aaa_1      &  (4\aaa_1)_I   &  3\aaa_1  \\
        &          &              &  2\aaa_1       &  \aaa_3  \\
        &          &              &                &  \aaa_1
\end{array}\right)$
 & $17$ &$2_6^{-2},4_3^{-3},3^{-1}$  \\
 & & $\subset 2\aaa_2\amalg 2\aaa_1\amalg \aaa_3$ & & \\
 \hline
& &
     $\left(\begin{array}{rrrrr}
2\aaa_1 & (4\aaa_1)_I  & (4\aaa_1)_I  &  (4\aaa_1)_{II}   &  3\aaa_1  \\
        & 2\aaa_1      & 2\aaa_2     &  (4\aaa_1)_I      &  3\aaa_1  \\
        &              & 2\aaa_1      &  (4\aaa_1)_I      &  3\aaa_1  \\
        &              &              &  2\aaa_1          &  \aaa_3  \\
        &              &              &                   &  \aaa_1
\end{array}\right)$
 & $17$ &$2_{II}^{+2},4_5^{-3},3^{-1}$  \\
 & & $\subset 2\aaa_1\amalg 2\aaa_2\amalg \aaa_3$ & & \\
 \hline
&  &
$(\aaa_1,2\aaa_1,2\aaa_1,2\aaa_1,2\aaa_1)\subset \aaa_1\amalg 2\aaa_2\amalg 2\aaa_2$  &
$17$ &$2_{3}^{+3},4_{II}^{+2},3^{+2}$     \\
\hline
& &
     $\left(\begin{array}{rrrrr}
\aaa_1 & \aaa_3   & 3\aaa_1         &  3\aaa_1      &  3\aaa_1  \\
       & 2\aaa_1  & (4\aaa_1)_{II}  &  (4\aaa_1)_I  &  (4\aaa_1)_I  \\
       &          & 2\aaa_1         & (4\aaa_1)_I   &  (4\aaa_1)_I  \\
       &          &                 & 2\aaa_1       & (4\aaa_1)_I \\
       &          &                 &               & 2\aaa_1
\end{array}\right)$
 & $17$ & $4_7^{+5}$     \\
 & & $\subset \aaa_3\amalg 6\aaa_1$ & & \\
\hline

\end{tabular}
\end{table}

\begin{table}
\begin{tabular}{|c||c|c|c|c|c|c|c|c|}
\hline
 {\bf n}& $G$  & $Deg$& $\rk S$ &$q_S$ \\
\hline
\hline
 $3$  & $(C_2)^2$ &
 $\left(\begin{array}{rrrrr}
2\aaa_1 & (4\aaa_1)_I   & (4\aaa_1)_I   & (4\aaa_1)_I      &  3\aaa_1 \\
        & 2\aaa_1       & 2\aaa_2       &  (4\aaa_1)_{II}  &  3\aaa_1  \\
        &               & \aaa_1        &  2\aaa_2         & 3\aaa_1      \\
        &               &               & 2\aaa_1          & \aaa_3 \\
       &            &                 &                    & \aaa_1
\end{array}\right)$
 & $17$ & $2_{II}^{-2},4_{3}^{-3}$     \\
 & & $\subset 2\aaa_1\amalg \aaa_7$ & & \\
 \hline
& &
     $\left(\begin{array}{rrrrr}
\aaa_1 & 3\aaa_1 & 3\aaa_1      &  3\aaa_1      &  5\aaa_1  \\
       & 2\aaa_1 & (4\aaa_1)_I  &  (4\aaa_1)_I  &  6\aaa_1  \\
       &         & 2\aaa_1      &  (4\aaa_1)_I  &  6\aaa_1  \\
       &         &              &  2\aaa_1      &  6\aaa_1 \\
       &         &          &            & 4\aaa_1
\end{array}\right)\subset 11\aaa_1$
 & $17$ & $2_1^{+1},4_5^{-3},8_1^{+1}$  \\
\hline
&  &
$(\aaa_1,2\aaa_1,2\aaa_1,2\aaa_1,4\aaa_1)\subset 3\aaa_1\amalg 2\aaa_2\amalg 4\aaa_1$  &
$17$ &$2_{3}^{-3},4_{5}^{-1},8_5^{-1},3^{-1}$     \\
\hline
& &
     $\left(\begin{array}{rrrrr}
\aaa_1 & \aaa_3   & 3\aaa_1        &  3\aaa_1           &  5\aaa_1  \\
       & 2\aaa_1  & (4\aaa_1)_{I}  &  (4\aaa_1)_{I}     &  6\aaa_1   \\
       &          & 2\aaa_1        &  (4\aaa_1)_I       &  6\aaa_1  \\
       &          &                &  2\aaa_1           &  6\aaa_1 \\
       &          &                &                    & 4\aaa_1
\end{array}\right)$
 & $17$ & $2_6^{+2},4_2^{-2},8_7^{+1}$   \\
& & $\subset \aaa_3\amalg 8\aaa_1$ &  &  \\
\hline

\end{tabular}
\end{table}

\begin{table}
\begin{tabular}{|c||c|c|c|c|c|c|c|c|}
\hline
 {\bf n}& $G$  & $Deg$& $\rk S$ &$q_S$ \\
\hline
\hline
 $3$  & $(C_2)^2$ &
     $\left(\begin{array}{rrrrr}
\aaa_1 & \aaa_3   & 3\aaa_1        &  5\aaa_1           &  3\aaa_1  \\
       & 2\aaa_1  & (4\aaa_1)_{I}  &  6\aaa_1           &  (4\aaa_1)_I   \\
       &          & 2\aaa_1        &  2\aaa_3           &  (4\aaa_1)_I  \\
       &          &                &  4\aaa_1           &  6\aaa_1 \\
       &          &                &                    & 2\aaa_1
\end{array}\right)$
 & $17$ & $2_{4}^{-2},4_7^{+3}$   \\
& & $\subset \aaa_3\amalg 2\aaa_3\amalg 2\aaa_1$ & & \\
\hline
& &
     $\left(\begin{array}{rrrrr}
\aaa_1 & \aaa_3   & 3\aaa_1        &  5\aaa_1           &  3\aaa_1  \\
       & 2\aaa_1  & (4\aaa_1)_{I}  &  6\aaa_1           &  (4\aaa_1)_{II}   \\
       &          & 2\aaa_1        &  2\aaa_3           &  (4\aaa_1)_I  \\
       &          &                &  4\aaa_1           &  6\aaa_1 \\
       &          &                &                    & 2\aaa_1
\end{array}\right)$
 & $17$ & $2_{II}^{+2},4_7^{+3}$   \\
& & $\subset \aaa_3\amalg 2\aaa_3\amalg 2\aaa_1$ & & \\
\hline
& &
     $\left(\begin{array}{rrrrr}
\aaa_1 & 3\aaa_1   & 3\aaa_1       &  3\aaa_1           &  5\aaa_1  \\
       & 2\aaa_1  & (4\aaa_1)_{I}  &  (4\aaa_1)_I       &  6\aaa_1   \\
       &          & 2\aaa_1        &  (4\aaa_1)_I       &  6\aaa_1  \\
       &          &                &  2\aaa_1           &  2\aaa_3 \\
       &          &                &                    &  4\aaa_1
\end{array}\right)\subset 5\aaa_1\amalg 2\aaa_3$
 & $17$ &$2_3^{-1},4_{0}^{+4}$  \\
\hline

\end{tabular}
\end{table}


\begin{table}
\begin{tabular}{|c||c|c|c|c|c|c|c|c|}
\hline
 {\bf n}& $G$  & $Deg$& $\rk S$ &$q_S$ \\
\hline
\hline
 $3$  & $(C_2)^2$ &
 $\left(\begin{array}{rrrrr}
4\aaa_1 & 6\aaa_1   & 6\aaa_1      &  6\aaa_1           &  5\aaa_1  \\
       & 2\aaa_1    & (4\aaa_1)_I  &  (4\aaa_1)_I       &  3\aaa_1   \\
       &            & 2\aaa_1      &  2\aaa_2           &  3\aaa_1  \\
       &          &                &  2\aaa_1           &  \aaa_3 \\
       &          &                &                    &  \aaa_1
\end{array}\right)$
 & $17$ & $2_{1}^{-3},4_3^{-1},8_3^{-1}$   \\
& & $\subset 6\aaa_1\amalg \aaa_5$ & & \\
\hline

&   &
$\left(\begin{array}{rrrrr}
4\aaa_1 & 2\aaa_3   & 6\aaa_1      &  6\aaa_1           &  5\aaa_1  \\
       & 2\aaa_1    & (4\aaa_1)_I  &  (4\aaa_1)_I       &  3\aaa_1   \\
       &            & 2\aaa_1      &  2\aaa_2           &  3\aaa_1  \\
       &          &                &  2\aaa_1           &  \aaa_3 \\
       &          &                &                    &  \aaa_1
\end{array}\right)$
 & $17$ & $2_{1}^{-3},4_2^{+2}$   \\
& & $\subset 2\aaa_3\amalg \aaa_5$ & & \\
\hline
&   &
$\left(\begin{array}{rrrrr}
4\aaa_1 & 6\aaa_1   & 6\aaa_1      &  6\aaa_1           &  5\aaa_1  \\
       & 2\aaa_1    & 2\aaa_2  &  (4\aaa_1)_{II}        &  3\aaa_1   \\
       &            & 2\aaa_1      &  2\aaa_2           &  3\aaa_1  \\
       &          &                &  2\aaa_1           &  \aaa_3 \\
       &          &                &                    &  \aaa_1
\end{array}\right)$
 & $17$ & $2_{II}^{+4},8_7^{+1}$   \\
& & $\subset 4\aaa_1\amalg \aaa_7$ & & \\
\hline
&  &
$(\aaa_1,(2\aaa_1,2\aaa_1)_I,4\aaa_1,4\aaa_1)\subset 13\aaa_1$  &
$17$ &$2_{5}^{-1},4_{6}^{-4}$     \\
\hline

\end{tabular}
\end{table}

\begin{table}
\begin{tabular}{|c||c|c|c|c|c|c|c|c|}
\hline
 {\bf n}& $G$  & $Deg$& $\rk S$ &$q_S$ \\
\hline
\hline
 $3$  & $(C_2)^2$ &
     $\left(\begin{array}{rrrrr}
2\aaa_1 & 3\aaa_1   & (4\aaa_1)_I      &  6\aaa_1           &  5\aaa_1  \\
        & \aaa_1    & \aaa_3           &  5\aaa_1           &  5\aaa_1  \\
        &           & 2\aaa_1          &  6\aaa_1           &  6\aaa_1  \\
        &          &                   &  4\aaa_1           &  4\aaa_2  \\
        &          &                &                       &  4\aaa_1
\end{array}\right)$
 & $17$ & $2_{6}^{+4},4_3^{-1},3^{+1}$   \\
& & $\subset 2\aaa_1\amalg \aaa_3\amalg 4\aaa_2$ & & \\
\hline
& &
     $\left(\begin{array}{rrrrr}
\aaa_1  & \aaa_3     & 5\aaa_1      &  3\aaa_1           &  5\aaa_1  \\
        & 2\aaa_1    & 6\aaa_1      &  (4\aaa_1)_I       &  6\aaa_1  \\
        &            & 4\aaa_1      &  6\aaa_1           &  8\aaa_1  \\
        &            &              &  2\aaa_1           &  2\aaa_3  \\
        &          &                &                    &  4\aaa_1
\end{array}\right)$
 & $17$ & $2_{6}^{+4},8_5^{-1}$   \\
& & $\subset \aaa_3\amalg 4\aaa_1\amalg 2\aaa_3$ & & \\
\hline
& &
     $\left(\begin{array}{rrrrr}
\aaa_1  & \aaa_3     & 3\aaa_1          &  5\aaa_1           &  5\aaa_1  \\
        & 2\aaa_1    & (4\aaa_1)_{II}   &  6\aaa_1           &  6\aaa_1  \\
        &            & 2\aaa_1          &  6\aaa_1           &  6\aaa_1  \\
        &            &                  &  4\aaa_1           &  8\aaa_1  \\
        &            &                  &                    &  4\aaa_1
\end{array}\right)$
 & $17$ & $2_{II}^{+2},4_7^{+3}$   \\
& & $\subset \aaa_3\amalg 2\aaa_1\amalg 8\aaa_1$ & & \\
\hline

\end{tabular}
\end{table}

\begin{table}
\begin{tabular}{|c||c|c|c|c|c|c|c|c|}
\hline
 {\bf n}& $G$  & $Deg$& $\rk S$ &$q_S$ \\
\hline
\hline
 $3$  & $(C_2)^2$ &
     $\left(\begin{array}{rrrrr}
4\aaa_1  & 4\aaa_2     & 6\aaa_1          &  6\aaa_1           &  5\aaa_1  \\
         & 4\aaa_1     & 2\aaa_3          &  6\aaa_1           &  5\aaa_1  \\
         &             & 2\aaa_1          &  (4\aaa_1)_I       &  3\aaa_1  \\
         &             &                  &  2\aaa_1           &  3\aaa_1  \\
        &            &                  &                      &  \aaa_1
\end{array}\right)$
 & $17$ & $2_{3}^{-3},4_{0}^{+2}$   \\
& & $\subset 2\aaa_5\amalg 3\aaa_1$    & & \\
\hline
&   &
     $\left(\begin{array}{rrrrr}
4\aaa_1  & 4\aaa_2     & 6\aaa_1          &  6\aaa_1           &  5\aaa_1  \\
         & 4\aaa_1     & 2\aaa_3          &  6\aaa_1           &  5\aaa_1  \\
         &             & 2\aaa_1          &  (4\aaa_1)_I       &  3\aaa_1  \\
         &             &                  &  2\aaa_1           &  \aaa_3  \\
        &            &                  &                      &  \aaa_1
\end{array}\right)$
 & $17$ & $2_{6}^{+4},4_{1}^{+1}$   \\
& & $\subset 2\aaa_5\amalg \aaa_3$    & & \\
\hline
&   &
     $\left(\begin{array}{rrrrr}
4\aaa_1  & 4\aaa_2     & 6\aaa_1          &  6\aaa_1           &  5\aaa_1  \\
         & 4\aaa_1     & 2\aaa_3          &  6\aaa_1           &  5\aaa_1  \\
         &             & 2\aaa_1          &  (4\aaa_1)_{II}    &  3\aaa_1  \\
         &             &                  &  2\aaa_1           &  \aaa_3  \\
        &            &                  &                      &  \aaa_1
\end{array}\right)$
 & $17$ & $2_{II}^{-4},4_{3}^{-1}$   \\
& & $\subset 2\aaa_5\amalg \aaa_3$    & & \\
\hline

\end{tabular}
\end{table}


\begin{table}
\begin{tabular}{|c||c|c|c|c|c|c|c|c|}
\hline
 {\bf n}& $G$  & $Deg$& $\rk S$ &$q_S$ \\
\hline
\hline
 $3$  & $(C_2)^2$ &
 $\left(\begin{array}{rrrrr}
2\aaa_1 & (4\aaa_1)_{II}   & (4\aaa_1)_I       &  (4\aaa_1)_I           &  (4\aaa_1)_I  \\
        & 2\aaa_1          & (4\aaa_1)_I       &  (4\aaa_1)_I           &  (4\aaa_1)_I   \\
       &                   & 2\aaa_1           &  (4\aaa_1)_I           &  (4\aaa_1)_I    \\
       &          &                            &  2\aaa_1               &  (4\aaa_1)_{II} \\
       &          &                            &                        &  2\aaa_1
\end{array}\right)$
 & $17$ &$4_7^{+5}$  \\
&  &  $\subset 10\aaa_1$  &  & \\
\hline
 &   &
 $\left(\begin{array}{rrrrr}
2\aaa_1 & (4\aaa_1)_{II}   & (4\aaa_1)_I       &  (4\aaa_1)_I           &  (4\aaa_1)_I  \\
        & 2\aaa_1          & (4\aaa_1)_I       &  (4\aaa_1)_I           &  (4\aaa_1)_I   \\
       &                   & 2\aaa_1           &  (4\aaa_1)_{II}        &  (4\aaa_1)_I    \\
       &          &                            &  2\aaa_1               &  (4\aaa_1)_{I} \\
       &          &                            &                        &  2\aaa_1
\end{array}\right)$
 & $17$ &$2_{II}^{+2},4_7^{+3}$  \\
&  &   $\subset 10\aaa_1$ & & \\
\hline
& &
     $\left(\begin{array}{rrrrr}
2\aaa_1 & (4\aaa_1)_I   & (4\aaa_1)_I      & (4\aaa_1)_I           &  (4\aaa_1)_I  \\
        & 2\aaa_1       & 2\aaa_2          & (4\aaa_1)_{II}        &  (4\aaa_1)_I  \\
        &               & 2\aaa_1          & 2\aaa_2               &  (4\aaa_1)_I  \\
        &               &                  & 2\aaa_1               &  (4\aaa_1)_I  \\
        &               &                  &                       &  2\aaa_1
\end{array}\right)$
 & $17$ & $2_{II}^{-2},4_2^{+2},8_5^{-1}$   \\
 &      & $\subset 2\aaa_1\amalg 2\aaa_3\amalg 2\aaa_1$ & & \\
\hline

\end{tabular}
\end{table}

\begin{table}
\begin{tabular}{|c||c|c|c|c|c|c|c|c|}
\hline
 {\bf n}& $G$  & $Deg$& $\rk S$ &$q_S$ \\
\hline
\hline
 $3$  & $(C_2)^2$ &
$(2\aaa_1,2\aaa_1,2\aaa_1,2\aaa_1,2\aaa_1)\subset 2\aaa_1\amalg 2\aaa_2\amalg 2\aaa_2$  &
$17$ &$4_{3}^{+3},3^{+2}$     \\
\hline
 &   &
 $\left(\begin{array}{rrrrr}
2\aaa_1 & 2\aaa_2          & (4\aaa_1)_{II}    &  (4\aaa_1)_I           &  (4\aaa_1)_I  \\
        & 2\aaa_1          & 2\aaa_2           &  (4\aaa_1)_I           &  (4\aaa_1)_I   \\
       &                   & 2\aaa_1           &  (4\aaa_1)_I           &  (4\aaa_1)_I    \\
       &          &                            &  2\aaa_1               &  (4\aaa_1)_{I} \\
       &          &                            &                        &  2\aaa_1
\end{array}\right)\subset 2\aaa_3\amalg 4\aaa_1$
 & $17$ &$2_{II}^{+2},4_0^{+2},8_7^{+1}$  \\
\hline
 &   &
 $\left(\begin{array}{rrrrr}
2\aaa_1 & 2\aaa_2          & (4\aaa_1)_{II}    &  (4\aaa_1)_I           &  (4\aaa_1)_I  \\
        & 2\aaa_1          & 2\aaa_2           &  (4\aaa_1)_I           &  (4\aaa_1)_I   \\
       &                   & 2\aaa_1           &  (4\aaa_1)_I           &  (4\aaa_1)_I    \\
       &          &                            &  2\aaa_1               &  (4\aaa_1)_{II} \\
       &          &                            &                        &  2\aaa_1
\end{array}\right)\subset 2\aaa_3\amalg 4\aaa_1$
 & $17$ &$2_{II}^{+4},8_7^{+1}$  \\
\hline
 &   &
 $\left(\begin{array}{rrrrr}
2\aaa_1 & 2\aaa_2          & (4\aaa_1)_{II}    &  (4\aaa_1)_I           &  (4\aaa_1)_I  \\
        & 2\aaa_1          & 2\aaa_2           &  (4\aaa_1)_I           &  (4\aaa_1)_I   \\
       &                   & 2\aaa_1           &  (4\aaa_1)_I           &  (4\aaa_1)_I   \\
       &          &                            &  2\aaa_1               &  2\aaa_2       \\
       &          &                            &                        &  2\aaa_1
\end{array}\right)\subset 2\aaa_3\amalg 2\aaa_2$
 & $17$ &$2_{II}^{-4},8_1^{+1},3^{-1}$  \\
\hline
\hline
\end{tabular}
\end{table}

\begin{table}
\begin{tabular}{|c||c|c|c|c|c|c|c|c|}
\hline
 {\bf n}& $G$  & $Deg$& $\rk S$ &$q_S$ \\
\hline
\hline
 $3$  & $(C_2)^2$ &
$\left(\begin{array}{rrrrr}
4\aaa_1 & 6\aaa_1          & 6\aaa_1           &  6\aaa_1               &  6\aaa_1  \\
        & 2\aaa_1          & 2\aaa_2           &  (4\aaa_1)_I           &  (4\aaa_1)_I   \\
       &                   & 2\aaa_1           &  (4\aaa_1)_I           &  (4\aaa_1)_I    \\
       &          &                            &  2\aaa_1               &  (4\aaa_1)_{II} \\
       &          &                            &                        &  2\aaa_1
\end{array}\right)$
 & $17$ &$2_{II}^{-4},8_1^{+1},3^{-1}$  \\
& & $\subset 4\aaa_1\amalg 2\aaa_2\amalg 4\aaa_1$ & & \\
\hline
& &
$\left(\begin{array}{rrrrr}
4\aaa_1 & 6\aaa_1          & 6\aaa_1           &  6\aaa_1               &  6\aaa_1        \\
        & 2\aaa_1          & 2\aaa_2           &  (4\aaa_1)_I           &  (4\aaa_1)_I    \\
       &                   & 2\aaa_1           &  (4\aaa_1)_I           &  (4\aaa_1)_I    \\
       &          &                            &  2\aaa_1               &  2\aaa_2        \\
       &          &                            &                        &  2\aaa_1
\end{array}\right)$
 & $17$ &$2_{II}^{-2},8_3^{-1},3^{+2}$  \\
& & $\subset 4\aaa_1\amalg 2\aaa_2\amalg 2\aaa_2$ & & \\
\hline
 &   &
 $\left(\begin{array}{rrrrr}
4\aaa_1 & 2\aaa_3          & 6\aaa_1           &  6\aaa_1               &  6\aaa_1  \\
        & 2\aaa_1          & (4\aaa_1)_{II}    &  (4\aaa_1)_I           &  (4\aaa_1)_I   \\
       &                   & 2\aaa_1           &  (4\aaa_1)_I           &  (4\aaa_1)_I    \\
       &          &                            &  2\aaa_1               &  (4\aaa_1)_{I} \\
       &          &                            &                        &  2\aaa_1
\end{array}\right)\subset 2\aaa_3\amalg 6\aaa_1$
 & $17$ &$2_{II}^{-2},4_3^{-3}$  \\
\hline

\end{tabular}
\end{table}


\begin{table}
\begin{tabular}{|c||c|c|c|c|c|c|c|c|}
\hline
 {\bf n}& $G$  & $Deg$& $\rk S$ &$q_S$ \\
\hline
\hline
 $3$  & $(C_2)^2$ &
  $\left(\begin{array}{rrrrr}
4\aaa_1 & 2\aaa_3          & 6\aaa_1           &  6\aaa_1               &  6\aaa_1        \\
        & 2\aaa_1          & (4\aaa_1)_I       &  (4\aaa_1)_I           &  (4\aaa_1)_{II}    \\
       &                   & 2\aaa_1           &  2\aaa_2               &  (4\aaa_1)_I    \\
       &          &                            &  2\aaa_1               &  (4\aaa_1)_I    \\
       &          &                            &                        &  2\aaa_1
\end{array}\right)$
 & $17$ &  $2_{II}^{+4},4_{5}^{-1},3^{-1}$  \\
& & $\subset 2\aaa_3\amalg 2\aaa_2\amalg 2\aaa_1$ & & \\
\hline
 &   &
 $\left(\begin{array}{rrrrr}
4\aaa_1 & 2\aaa_3          & 6\aaa_1           &  6\aaa_1               &  6\aaa_1  \\
        & 2\aaa_1          & (4\aaa_1)_{I}    &  (4\aaa_1)_I           &  (4\aaa_1)_I   \\
       &                   & 2\aaa_1           &  2\aaa_2              &  (4\aaa_1)_{II}    \\
       &          &                            &  2\aaa_1              &  2\aaa_2 \\
       &          &                            &                       &  2\aaa_1
\end{array}\right)$
 & $17$ &$2_{II}^{+4},8_7^{+1}$  \\
 &  &  $\subset 2\aaa_3\amalg 2\aaa_3$ &  & \\
 \hline
  &   &
 $\left(\begin{array}{rrrrr}
4\aaa_1 & 6\aaa_1          & 6\aaa_1           &  6\aaa_1               &  6\aaa_1  \\
        & 2\aaa_1          & (4\aaa_1)_{I}    &  (4\aaa_1)_I           &  (4\aaa_1)_I   \\
       &                   & 2\aaa_1           &  2\aaa_2           &  (4\aaa_1)_{II}    \\
       &          &                            &  2\aaa_1               &  2\aaa_2 \\
       &          &                            &                        &  2\aaa_1
\end{array}\right)$
 & $17$ &$2_{II}^{+2},4_7^{+3}$  \\
 &  &  $\subset 6\aaa_1\amalg 2\aaa_3$ &  & \\
\hline

\end{tabular}
\end{table}

\begin{table}
\begin{tabular}{|c||c|c|c|c|c|c|c|c|}
\hline
 {\bf n}& $G$  & $Deg$& $\rk S$ &$q_S$ \\
\hline
\hline
 $3$  & $(C_2)^2$ &
 $\left(\begin{array}{rrrrr}
4\aaa_1 & 2\aaa_3          & 6\aaa_1           &  6\aaa_1               &  6\aaa_1  \\
        & 2\aaa_1          & 2\aaa_2           &  (4\aaa_1)_{II}        &  (4\aaa_1)_I   \\
        &                  & 2\aaa_1           &  2\aaa_2               &  (4\aaa_1)_{I}    \\
       &          &                            &  2\aaa_1               &  (4\aaa_1)_I \\
       &          &                            &                        &  2\aaa_1
\end{array}\right)\subset 2\ddd_5\amalg 2\aaa_1$
 & $17$ &$2_{II}^{+4},4_7^{+1}$  \\
\hline
&   &
 $\left(\begin{array}{rrrrr}
2\aaa_1 & (4\aaa_1)_I        & (4\aaa_1)_I           &  6\aaa_1       &  6\aaa_1  \\
        & 2\aaa_1            & (4\aaa_1)_{I}         &  6\aaa_1       &  6\aaa_1  \\
       &                     & 2\aaa_1               &  6\aaa_1       &  6\aaa_1  \\
       &                     &                       &  4\aaa_1       &  8\aaa_1  \\
       &                     &                       &                &  4\aaa_1
\end{array}\right)\subset 14\aaa_1$
 & $17$ &$2_{4}^{-2},4_7^{-3}$  \\
\hline
 &   &
 $\left(\begin{array}{rrrrr}
2\aaa_1 & (4\aaa_1)_{II}        & (4\aaa_1)_I           &  6\aaa_1       &  6\aaa_1  \\
        & 2\aaa_1            & (4\aaa_1)_{I}         &  6\aaa_1       &  6\aaa_1  \\
       &                     & 2\aaa_1               &  6\aaa_1       &  6\aaa_1  \\
       &                     &                       &  4\aaa_1       &  8\aaa_1  \\
       &                     &                       &                &  4\aaa_1
\end{array}\right)\subset 14\aaa_1$
 & $17$ &$2_{II}^{+2},4_7^{+3}$  \\
\hline

\end{tabular}
\end{table}

\begin{table}
\begin{tabular}{|c||c|c|c|c|c|c|c|c|}
\hline
 {\bf n}& $G$  & $Deg$& $\rk S$ &$q_S$ \\
\hline
\hline
 $3$  & $(C_2)^2$ &
$\left(\begin{array}{rrrrr}
2\aaa_1 & (4\aaa_1)_I          & 6\aaa_1       &  (4\aaa_1)_I          &  6\aaa_1  \\
        & 2\aaa_1              & 2\aaa_3       &  (4\aaa_1)_{II}       &  6\aaa_1  \\
        &                      & 4\aaa_1       &  6\aaa_1              &  8\aaa_1  \\
        &          &                            &  2\aaa_1             &  2\aaa_3 \\
        &          &                            &                      &  4\aaa_1
\end{array}\right)$
 & $17$ & $2_{II}^{-4},4_{3}^{-1}$   \\
 &  &  $\subset 2\aaa_1\amalg 2\aaa_3\amalg 2\aaa_3$ & & \\
 \hline
 &  &
 $\left(\begin{array}{rrrrr}
4\aaa_1 & 2\aaa_3          & 6\aaa_1       &  8\aaa_1       &  6\aaa_1  \\
        & 2\aaa_1          & 2\aaa_2       &  6\aaa_1       &  (4\aaa_1)_I  \\
        &                  & 2\aaa_1       &  6\aaa_1       &  (4\aaa_1)_I  \\
        &          &                       &  4\aaa_1       &  6\aaa_1 \\
        &          &                            &           &  2\aaa_1
\end{array}\right)$
 & $17$ & $4_{7}^{+3}$   \\
 &  &  $\subset 2\ddd_4\amalg 6\aaa_1$ & & \\
 \hline
   &  &
 $\left(\begin{array}{rrrrr}
4\aaa_1 & 2\aaa_3          & 6\aaa_1       &  8\aaa_1       &  6\aaa_1  \\
        & 2\aaa_1          & 2\aaa_2       &  6\aaa_1       &  (4\aaa_1)_I  \\
        &                  & 2\aaa_1       &  6\aaa_1       &  (4\aaa_1)_I  \\
        &          &                       &  4\aaa_1       &  2\aaa_3 \\
        &          &                            &           &  2\aaa_1
\end{array}\right)$
 & $17$ & $2_{II}^{+2},8_{7}^{+1}$   \\
 &  &  $\subset 2\ddd_4\amalg 2\aaa_3$ & & \\
\hline

\end{tabular}
\end{table}

\begin{table}
\begin{tabular}{|c||c|c|c|c|c|c|c|c|}
\hline
 {\bf n}& $G$  & $Deg$& $\rk S$ &$q_S$ \\
\hline
\hline
 $3$  & $(C_2)^2$ &
 $\left(\begin{array}{rrrrr}
4\aaa_1 & 4\aaa_2          & 6\aaa_1       &  6\aaa_1       &  6\aaa_1  \\
        & 4\aaa_1          & 2\aaa_3       &  6\aaa_1       &  6\aaa_1  \\
        &                  & 2\aaa_1       &  2\aaa_2       &  (4\aaa_1)_{I}  \\
        &          &                       &  2\aaa_1       &  (4\aaa_1)_I \\
        &          &                            &           &  2\aaa_1
\end{array}\right)$
 & $17$ & $2_{II}^{+2}, 4_7^{+1}$   \\
 &  &  $\subset 2\eee_6\amalg 2\aaa_1$ & & \\
 \hline
&  &
 $\left(\begin{array}{rrrrr}
4\aaa_1 & 4\aaa_2          & 8\aaa_1       &  6\aaa_1       &  6\aaa_1  \\
        & 4\aaa_1          & 4\aaa_2       &  6\aaa_1       &  6\aaa_1  \\
        &                  & 4\aaa_1       &  2\aaa_3       &  6\aaa_1  \\
        &          &                       &  2\aaa_1       &  (4\aaa_1)_I \\
        &          &                            &           &  2\aaa_1
\end{array}\right)$
 & $17$ & $2_{II}^{-2}, 4_3^{-1}$   \\
 &  &  $\subset 2\aaa_7\amalg 2\aaa_1$ & & \\
 \hline
&   &
$\left(\begin{array}{rrrrr}
4\aaa_1 & 2\aaa_3               & 8\aaa_1      &  6\aaa_1              &  8\aaa_1  \\
        & 2\aaa_1              & 6\aaa_1       &  (4\aaa_1)_{I}        &  6\aaa_1  \\
        &                      & 4\aaa_1       &  2\aaa_3              &  8\aaa_1  \\
        &                      &               &  2\aaa_1              &  6\aaa_1 \\
        &                      &               &                       &  4\aaa_1
\end{array}\right)$
 & $17$ & $2_{II}^{+2},8_{7}^{+1}$   \\
 &  &  $\subset 2\aaa_3\amalg 2\aaa_3\amalg 4\aaa_1$ & & \\
 \hline
  &    &
$((2\aaa_1,2\aaa_1)_I,4\aaa_1,4\aaa_1,4\aaa_1)\subset 4\aaa_1\amalg 4\aaa_3$  &
$17$ &   $4_{7}^{+3}$     \\
\hline

\end{tabular}
\end{table}


\begin{table}
\begin{tabular}{|c||c|c|c|c|c|c|c|c|}
\hline
 {\bf n}& $G$  & $Deg$& $\rk S$ &$q_S$ \\
\hline
\hline
 $3$  & $(C_2)^2$ &
 $\left(\begin{array}{rrrrrr}
2\aaa_1   & (4\aaa_1)_{I}    & (4\aaa_1)_I    &  3\aaa_1       &  (4\aaa_1)_I   & 3\aaa_1     \\
          & 2\aaa_1          & (4\aaa_1)_I    &  3\aaa_1       &  (4\aaa_1)_I   & 3\aaa_1     \\
          &                  & 2\aaa_1        &  \aaa_3        &  (4\aaa_1)_{II}  & 3\aaa_1     \\
          &                  &                &  \aaa_1        &  3\aaa_1       & 2\aaa_1     \\
          &                  &                &                &  2\aaa_1       & \aaa_3       \\
          &                  &                &                &                & \aaa_1
\end{array}\right)$
 & $18$ &$4_6^{+4}$  \\
&  &
$\subset 4\aaa_1\amalg \aaa_3\amalg \aaa_3$ &  &  \\
\hline
&   &
$\left(\begin{array}{rrrrrr}
2\aaa_1   & 2\aaa_2          & 3\aaa_1        &  (4\aaa_1)_I   &  (4\aaa_1)_I   & 3\aaa_1     \\
          & 2\aaa_1          & \aaa_3         &  (4\aaa_1)_I   &  (4\aaa_1)_I   & 3\aaa_1     \\
          &                  & \aaa_1         &  3\aaa_1       &  3\aaa_1       & 2\aaa_1     \\
          &                  &                &  2\aaa_1       &  (4\aaa_1)_I   & 3\aaa_1     \\
          &                  &                &                &  2\aaa_1       & 3\aaa_1     \\
          &                  &                &                &                & \aaa_1
\end{array}\right)$
 & $18$ &$4_6^{+4}$  \\
&  &
$\subset \aaa_5\amalg 5\aaa_1$ &  &  \\
\hline

\end{tabular}
\end{table}

\begin{table}
\begin{tabular}{|c||c|c|c|c|c|c|c|c|}
\hline
 {\bf n}& $G$  & $Deg$& $\rk S$ &$q_S$ \\
\hline
\hline
 $3$  & $(C_2)^2$ &
$\left(\begin{array}{rrrrrr}
\aaa_1    & 2\aaa_1          & 3\aaa_1        &  3\aaa_1       &  3\aaa_1       & 5\aaa_1     \\
          & \aaa_1           & \aaa_3         &  3\aaa_1       &  3\aaa_1       & 5\aaa_1     \\
          &                  & 2\aaa_1        &  (4\aaa_1)_I   &  (4\aaa_1)_{I} & 6\aaa_1     \\
          &                  &                &  2\aaa_1       &  (4\aaa_1)_I   & 6\aaa_1     \\
          &                  &                &                &  2\aaa_1       & 2\aaa_3       \\
          &                  &                &                &                & 4\aaa_1
\end{array}\right)$
 & $18$ &$2_7^{+1},4_7^{+3}$  \\
 &  &  $\subset \aaa_1\amalg \aaa_3\amalg 2\aaa_1\amalg 2\aaa_3$  & & \\
\hline
&   &
$\left(\begin{array}{rrrrrr}
\aaa_1    & 3\aaa_1          & 3\aaa_1        &  3\aaa_1       &  2\aaa_1       & 5\aaa_1     \\
          & 2\aaa_1          & (4\aaa_1)_I    &  (4\aaa_1)_I   &  3\aaa_1       & 6\aaa_1     \\
          &                  & 2\aaa_1        &  (4\aaa_1)_I   &  3\aaa_1       & 6\aaa_1     \\
          &                  &                &  2\aaa_1       &  \aaa_3        & 6\aaa_1     \\
          &                  &                &                &  \aaa_1        & 5\aaa_1      \\
          &                  &                &                &                & 4\aaa_1
\end{array}\right)$
 & $18$ &$2_5^{-1},4_6^{+2},8_7^{+1}$  \\
 &  &  $\subset 5\aaa_1\amalg \aaa_3\amalg 4\aaa_1$  & & \\
\hline

\end{tabular}
\end{table}


\begin{table}
\begin{tabular}{|c||c|c|c|c|c|c|c|c|}
\hline
 {\bf n}& $G$  & $Deg$ & $\rk S$ &$q_S$ \\
\hline
\hline
 $3$  & $(C_2)^2$ &
 $\left(\begin{array}{rrrrrr}
4\aaa_1   & 2\aaa_3          & 6\aaa_1        &  5\aaa_1       &  6\aaa_1         & 5\aaa_1     \\
          & 2\aaa_1          & (4\aaa_1)_I    &  3\aaa_1       &  (4\aaa_1)_I     & 3\aaa_1     \\
          &                  & 2\aaa_1        &  \aaa_3        &  (4\aaa_1)_{II}  & 3\aaa_1     \\
          &                  &                &  \aaa_1        &  3\aaa_1         & 2\aaa_1     \\
          &                  &                &                &  2\aaa_1         & \aaa_3       \\
          &                  &                &                &                  & \aaa_1
\end{array}\right)$
 & $18$ &$2_{II}^{-2},4_2^{-2}$  \\
&  &
$\subset 2\aaa_3\amalg \aaa_3\amalg \aaa_3$ &  &  \\
\hline
& &
 $\left(\begin{array}{rrrrrr}
4\aaa_1   & 6\aaa_1          & 6\aaa_1        &  6\aaa_1       &  5\aaa_1         & 5\aaa_1     \\
          & 2\aaa_1          & (4\aaa_1)_I    &  (4\aaa_1)_I   &  3\aaa_1         & 3\aaa_1     \\
          &                  & 2\aaa_1        &  2\aaa_2       &  3\aaa_1         & 3\aaa_1     \\
          &                  &                &  2\aaa_1       &  \aaa_3          & 3\aaa_1     \\
          &                  &                &                &  \aaa_1          & 2\aaa_1       \\
          &                  &                &                &                  & \aaa_1
\end{array}\right)$
 & $18$ &$2_{2}^{-2},4_7^{+1},8_1^{+1}$  \\
&  &
$\subset 6\aaa_1\amalg \aaa_5\amalg \aaa_1$ &  &  \\
\hline

\end{tabular}
\end{table}

\begin{table}
\begin{tabular}{|c||c|c|c|c|c|c|c|c|}
\hline
 {\bf n}& $G$  & $Deg$ & $\rk S$ &$q_S$ \\
\hline
\hline
 $3$  & $(C_2)^2$ &

$(\aaa_1,\aaa_1, (2\aaa_1,2\aaa_1)_I, 4\aaa_1,4\aaa_1)\subset 6\aaa_1\amalg 4\aaa_2$
& $18$ &$2_2^{+2},4_6^{-2},3^{+1}$ \\
\hline
&  &
 $\left(\begin{array}{rrrrrr}
\aaa_1    & \aaa_3           & 2\aaa_1         &  3\aaa_1       &  5\aaa_1         & 5\aaa_1     \\
          & 2\aaa_1          & 3\aaa_1         &  (4\aaa_1)_I   &  6\aaa_1         & 6\aaa_1     \\
          &                  & \aaa_1          &  \aaa_3        &  5\aaa_1         & 5\aaa_1     \\
          &                  &                 & 2\aaa_1        &  6\aaa_1         & 6\aaa_1     \\
          &                  &                 &                &  4\aaa_1         & 8\aaa_1     \\
          &                  &                 &                &                  & 4\aaa_1
\end{array}\right)$
 & $18$ &$2_{6}^{-2},4_4^{-2}$  \\
&  &
$\subset \aaa_3\amalg \aaa_3\amalg 8\aaa_1$ &  &  \\
\hline
&  &
 $\left(\begin{array}{rrrrrr}
\aaa_1    & \aaa_3           & 2\aaa_1         &  3\aaa_1       &  5\aaa_1         & 5\aaa_1     \\
          & 2\aaa_1          & 3\aaa_1         &  (4\aaa_1)_{II}&  6\aaa_1         & 6\aaa_1     \\
          &                  & \aaa_1          &  \aaa_3        &  5\aaa_1         & 5\aaa_1     \\
          &                  &                 & 2\aaa_1        &  6\aaa_1         & 6\aaa_1     \\
          &                  &                 &                &  4\aaa_1         & 8\aaa_1     \\
          &                  &                 &                &                  & 4\aaa_1
\end{array}\right)$
 & $18$ &$2_{II}^{-2},4_2^{-2}$  \\
&  &
$\subset \aaa_3\amalg \aaa_3\amalg 8\aaa_1$ &  &  \\
\hline

\end{tabular}%
\end{table}

\vskip1cm

\newpage

\section{Final remarks }
\label{sec5.remarks}

We hope to consider similar classification of degenerations and Picard
lattices $S$ of K\"ahlerian K3 surfaces with
remaining small symplectic automorphism groups $(C_2)^2$, $C_3$, $C_2$ and $C_1$
later as well.
Now, for these groups we obtain

\begin{theorem}
If the Picard lattice $S=MS_X$ of a K3 surface $X$ with
finite symplectic automorphism group and with at least one $-2$ curve
is different from all lattices of lines of Tables 1---4 of Section 4
from \cite{Nik12}, lattices of Table 1 from \cite{Nik13} and lattices
of Table 1 of this paper which are not denoted by $o$ (for example, if its genus
is different), then the symplectic automorphism group of $X$ is one of groups
$(C_2)^2$, $C_3$, $C_2$ or $C_1$.
\end{theorem}


\vfill



\begin{table}
\label{table4}
\caption{Types and lattices $S$ of degenerations of
codimension $1$ of K\"ahlerian K3 surfaces
with finite symplectic automorphism groups $G=Clos(G)$ (from \cite{Nik10}).}



\begin{tabular}{|c||c|c|c|c|c|c|c|c|}
\hline
 {\bf n}& $|G|$& $i$&  $G$  & $\rk S_G$ &$q_{S_G}$  &$Deg$& $\rk S$ &$q_S$ \\
\hline
\hline
 $1$  &$2$ & $1$& $C_2$ &   $8$     &$2_{II}^{+8}$ & $\aaa_1$   &  $9$  &
$2_7^{+9}$\\
\hline
      &    &    &    &              &              & $2\aaa_1$  &  $9$  &
$2_{II}^{-6},4_3^{-1}$\\
\hline
\hline
 $2$  & $3$& $1$&$C_3$  &    $12$   & $3^{+6}$       & $\aaa_1$   &  $13$ &
$2_3^{-1},3^{+6}$\\
\hline
      &    &    &       &           &                & $3\aaa_1$  & $13$  &
$2_1^{+1},3^{-5}$\\
\hline
\hline
 $3$  &$4$ & $2$&$C_2^2$&    $12$   & $2_{II}^{-6},4_{II}^{-2}$& $\aaa_1$   & $13$  &
$2_3^{+7},4_{II}^{+2}$\\
\hline
      &    &    &       &           &                            & $2\aaa_1$& $13$  &
$2_{II}^{-4},4_7^{-3}$\\
\hline
      &    &    &       &           &                & $4\aaa_1$  & $13$  &
$2_{II}^{-6},8_3^{-1}$\\
\hline
\hline
$4$   &$4$& $1$ &$C_4$  &   $14$     &$2_2^{+2},4_{II}^{+4}$ & $\aaa_1$   & $15$  &
$2_5^{-3},4_{II}^{+4}$\\
\hline
      &   &     &       &           &                        & $2\aaa_1$  & $15$  &
$4_1^{-5}$\\
\hline
      &   &     &       &           &                        & $4\aaa_1$  & $15$  &
$2_2^{+2},4_{II}^{+2},8_7^{+1}$\\
\hline
      &   &     &       &           &                        & $\aaa_2$   & $15$  &
$2_1^{+1},4_{II}^{-4}$\\
\hline
\hline
 $6$  &$6$& $1$ & $D_6$ & $14$      &$2_{II}^{-2},3^{+5}$    &  $\aaa_1$   & $15$  &
$2_7^{-3},3^{+5}$\\
\hline
      &   &     &       &           &                        &  $2\aaa_1$  & $15$  &
$4_3^{-1},3^{+5}$\\
\hline
      &   &     &       &           &                        &  $3\aaa_1$  & $15$  &
$2_1^{-3},3^{-4}$\\
\hline
      &   &     &       &           &                        &  $6\aaa_1$  & $15$  &
$4_1^{+1},3^{+4}$\\
\hline
\hline
  $9$ &$8$& $5$ &$C_2^3$&  $14$     & $2_{II}^{+6},4_2^{+2}$ &  $2\aaa_1$  & $15$  &
$2_{II}^{-4},4_5^{-3}$\\
\hline
      &   &     &       &           &               &  $4\aaa_1$  & $15$  &
$2_{II}^{+6},8_1^{+1}$\\
\hline
      &   &     &       &           &               &  $8\aaa_1$  & $15$  &
$2_{II}^{+6},4_1^{+1}$\\
\hline
\hline
  $10$&$8$& $3$ &$D_8$  &  $15$     & $4_1^{+5}$    &  $\aaa_1$  & $16$   &
$2_1^{+1},4_7^{+5}$\\
\hline
      &   &     &       &           &               &$(2\aaa_1)_I$& $16$   &
$2_6^{-2},4_6^{-4}$\\
\hline
      &   &     &       &           &               &$(2\aaa_1)_{II}$& $16$   &
$2_{II}^{+2},4_{II}^{+4}$\\
\hline
      &   &     &       &           &               &  $4\aaa_1$ & $16$  &
$4_7^{+3},8_1^{+1}$\\
\hline
      &   &     &       &           &               &  $8\aaa_1$ & $16$  &
$4_0^{+4}$\\
\hline
      &   &     &       &           &               &  $2\aaa_2$ & $16$  &
$4_{II}^{+4}$\\
\hline
\hline
  $12$&$8$& $4$ & $Q_8$ &  $17$     & $2_7^{-3},8_{II}^{-2}$ &  $8\aaa_1$ & $18$  &
$2_7^{-3},16_3^{-1}$\\
\hline
      &   &     &       &           &                        &  $\aaa_2$ & $18$   &
$2_6^{-2},8_{II}^{-2}$\\
\hline
\hline
 $16$ &$10$&$1$ &$D_{10}$&  $16$    & $5^{+4}$      &  $\aaa_1$  & $17$  &
$2_7^{+1},5^{+4}$\\
\hline
      &   &     &       &           &               &  $5\aaa_1$  &$17$  &
$2_7^{+1},5^{-3}$\\
\hline
\end{tabular}
\end{table}

\begin{table}

\begin{tabular}{|c||c|c|c|c|c|c|c|c|}
\hline
 {\bf n}& $|G|$& $i$&  $G$   & $\rk S_G$ &$q_{S_G}$ &$Deg$& $\rk S$ &$q_S$ \\
\hline
\hline
 $17$&$12$ & $3$&$\AAA_4$&  $16$& $2_{II}^{-2},4_{II}^{-2},3^{+2}$ & $\aaa_1$ &$17$  &
$2_7^{-3},4_{II}^{+2},3^{+2}$\\
\hline
     &     &    &        &           &     &$3\aaa_1$ &  $17$  &
$2_1^{-3},4_{II}^{+2},3^{-1}$\\
\hline
     &     &    &        &           &            &$4\aaa_1$ &  $17$  &
$2_{II}^{-2},8_3^{-1},3^{+2}$\\
\hline
     &     &    &        &           &             &$6\aaa_1$ &  $17$  &
$4_1^{-3},3^{+1}$\\
\hline
     &     &    &        &           &             &$12\aaa_1$&  $17$  &
$2_{II}^{-2},8_1^{+1},3^{-1}$\\
\hline
\hline
 $18$&$12$ &$4$ &$D_{12}$&  $16$&$2_{II}^{+4},3^{+4}$ &  $\aaa_1$&  $17$&
$2_7^{+5},3^{+4}$\\
\hline
     &     &    &        &           &   &  $2\aaa_1$    &  $17$  &
$2_{II}^{+2},4_7^{+1},3^{+4}$\\
\hline
     &     &    &        &           &   &  $3\aaa_1$    &  $17$  &
$2_5^{-5},3^{-3}$\\
\hline
     &     &    &        &           &   &  $6\aaa_1$    &  $17$  &
$2_{II}^{-2},4_1^{+1},3^{+3}$\\
\hline
\hline
 $21$&$16$ &$14$&$C_2^4$ &  $15$ & $2_{II}^{+6},8_I^{+1}$ &  $4\aaa_1$ &  $16$  &
$2_{II}^{+4},4_{II}^{+2}$\\
\hline
     &     &    &        &           &          &$16\aaa_1$&  $16$  &
$2_{II}^{+6}$\\
\hline
\hline
 $22$& $16$&$11$&$C_2\times D_8$& $16$& $2_{II}^{+2},4_0^{+4}$ &$2\aaa_1$&$17$  &
$4_7^{+5}$\\
\hline
     &     &    &        &           &     &  $4\aaa_1$    &  $17$  &
$2_{II}^{+2},4_0^{+2},8_7^{+1}$\\
\hline
     &     &    &        &           &     &  $8\aaa_1$    &  $17$  &
$2_{II}^{+2},4_7^{+3}$\\
\hline
\hline
 $26$& $16$& $8$&$SD_{16}$& $18$ &$2_7^{+1},4_7^{+1},8_{II}^{+2}$&$8\aaa_1$&$19$  &
$2_1^{+1},4_1^{+1},16_3^{-1}$\\
\hline
     &     &    &        &       &              &  $2\aaa_2$    &  $19$  &
$2_5^{-1},8_{II}^{-2}$\\
\hline
\hline
$30$ & $18$& $4$&$\AAA_{3,3}$&$16$   &$3^{+4},9^{-1}$&  $3\aaa_1$ &  $17$  &
$2_5^{-1},3^{-3},9^{-1}$\\
\hline
     &     &    &        &           &               &  $9\aaa_1$ &  $17$  &
$2_3^{-1},3^{+4}$\\
\hline
\hline
 $32$& $20$& $3$&$Hol(C_5)$& $18$    &$2_6^{-2},5^{+3}$&  $2\aaa_1$ &$19$ &
$4_1^{+1},5^{+3}$\\
\hline
     &     &    &        &           &                 &  $5\aaa_1$  &$19$&
$2_1^{+3},5^{-2}$\\
\hline
     &     &    &        &           &                 &  $10\aaa_1$ &$19$  &
$4_5^{-1},5^{+2}$\\
\hline
     &     &    &        &           &               &  $5\aaa_2$   &$19$  &
$2_5^{-1},5^{-2}$\\
\hline
\hline
$33$ & $21$&$1$ &$C_7\rtimes C_3$&$18$&  $7^{+3}$    &  $7\aaa_1$    &  $19$  &
$2_1^{+1},7^{+2}$\\
\hline
\hline
 $34$& $24$&$12$&$\SSS_4$&  $17$      & $4_3^{+3},3^{+2}$& $\aaa_1$&  $18$  &
$2_5^{-1},4_1^{+3},3^{+2}$\\
\hline
     &     &    &        &            &                  & $2\aaa_1$    &  $18$  &
$2_2^{+2},4_{II}^{+2},3^{+2}$\\
\hline
     &     &    &        &           &                   & $3\aaa_1$    &  $18$  &
$2_7^{+1},4_5^{-3},3^{-1}$ \\
\hline
     &     &    &        &           &                   & $4\aaa_1$    &  $18$  &
$4_3^{-1},8_3^{-1},3^{+2}$ \\
\hline
     &     &    &        &           &                   & $(6\aaa_1)_I$&  $18$  &
$2_4^{-2},4_0^{+2},3^{+1}$\\
\hline
     &     &    &        &           &                   & $(6\aaa_1)_{II}$&$18$ &
$2_{II}^{+2},4_{II}^{-2},3^{+1}$ \\
\hline
     &     &    &        &           &                   & $8\aaa_1$       &$18$  &
$4_2^{+2},3^{+2}$\\
\hline
     &     &    &        &           &                   & $12\aaa_1$    & $18$  &
$4_5^{-1},8_7^{+1},3^{-1}$\\
\hline
     &     &    &        &           &                   & $6\aaa_2$    &  $18$  &
$4_{II}^{-2},3^{+1}$\\
\hline

\end{tabular}

\end{table}

\begin{table}

\begin{tabular}{|c||c|c|c|c|c|c|c|c|}
\hline
 {\bf n}& $|G|$& $i$&  $G$   & $\rk S_G$ &$q_{S_G}$       &$Deg$& $\rk S$ &$q_S$    \\
\hline
\hline
 $39$&$32$ &$27$&$2^4C_2$&  $17$     &$2_{II}^{+2},4_0^{+2},8_7^{+1}$ & $4\aaa_1$&$18$&
$4_6^{+4}$\\
\hline
     &     &    &        &           &                                & $8\aaa_1$&  $18$&
$2_{II}^{+2},4_7^{+1},8_7^{+1}$\\
\hline
     &     &    &        &           &                                & $16\aaa_1$&$18$&
$2_{II}^{+2},4_6^{+2}$\\
\hline
\hline
$40$ &$32$ &$49$&$Q_8*Q_8$& $17$     & $4_7^{+5}$ & $8\aaa_1$  & $18$ &
$4_6^{+4}$\\
\hline
\hline
$46$ &$36$ &$9$ &$3^2C_4$ & $18$     & $2_6^{-2},3^{+2},9^{-1}$& $6\aaa_1$   &  $19$  &
$4_7^{+1},3^{+1},9^{-1}$               \\
\hline
     &     &    &        &           &                           & $9\aaa_1$ &  $19$  &
$2_5^{-3},3^{+2}$ \\
\hline
     &     &    &        &           &                           & $9\aaa_2$ &  $19$  &
$2_5^{-1},3^{+2}$\\
\hline
\hline
$48$ &$36$ &$10$&$\SSS_{3,3}$&$18$   & $2_{II}^{-2},3^{+3},9^{-1}$ & $3\aaa_1$& $19$  &
$2_5^{+3},3^{-2},9^{-1}$\\
\hline
     &     &    &        &           &                             & $6\aaa_1$& $19$  &
$4_1^{+1},3^{+2},9^{-1}$\\
\hline
     &     &    &        &           &                             & $9\aaa_1$& $19$  &
$2_7^{-3},3^{+3}$\\
\hline
\hline
$49$ & $48$&$50$&$2^4C_3$& $17$      & $2_{II}^{-4},8_1^{+1},3^{-1}$&$4\aaa_1$&$18$  &
$2_{II}^{-2},4_{II}^{+2},3^{-1}$\\
\hline
     &     &    &        &           &                              &$12\aaa_1$&$18$  &
$2_{II}^{-2},4_2^{-2}$\\
\hline
     &     &    &        &           &                              &$16\aaa_1$&$18$  &
$2_{II}^{-4},3^{-1}$ \\
\hline
\hline
$51$ & $48$&$48$&$C_2\times \SSS_4$&$18$&$2_{II}^{+2},4_2^{+2},3^{+2}$ & $2\aaa_1$ &  $19$  &
$4_1^{+3},3^{+2}$\\
\hline
     &     &    &        &           &                               & $4\aaa_1$   &  $19$  &
$2_{II}^{+2},8_1^{+1},3^{+2}$\\
\hline
     &     &    &        &           &                               & $6\aaa_1$   &  $19$  &
$4_7^{-3},3^{+1}$\\
\hline
     &     &    &        &           &                               & $8\aaa_1$   &  $19$  &
$2_{II}^{-2},4_5^{-1},3^{+2}$\\
\hline
     &     &    &        &           &                               &$12\aaa_1$   &  $19$  &
$2_{II}^{-2},8_7^{+1},3^{-1}$\\
\hline
\hline
$55$ & $60$&$5$ &$\AAA_5$ & $18$      &$2_{II}^{-2},3^{+1},5^{-2}$& $\aaa_1$& $19$&
$2_7^{-3},3^{+1},5^{-2}$               \\
\hline
     &     &    &        &           &                            &$5\aaa_1$& $19$&
$2_3^{+3},3^{+1},5^{+1}$               \\
\hline
     &     &    &        &           &                            &$6\aaa_1$& $19$&
$4_1^{+1},5^{-2}$\\
\hline
     &     &    &        &           &                            &$10\aaa_1$&$19$&
$4_7^{+1},3^{+1},5^{-1}$\\
\hline
     &     &    &        &           &                            &$15\aaa_1$&  $19$&
$2_5^{+3},5^{-1}$\\
\hline
\hline
$56$ &$64$&$138$&$\Gamma_{25}a_1$&$18$& $4_5^{+3},8_1^{+1}$ & $8\aaa_1$ &  $19$  &
$4_4^{-2},8_5^{-1}$               \\
\hline
     &     &    &        &           &                      &$16\aaa_1$ &  $19$  &
$4_5^{+3}$ \\
\hline
\hline
$61$ &$72$ &$43$&$\AAA_{4,3}$& $18$  & $4_{II}^{-2},3^{-3}$&$3\aaa_1$&$19$&
$2_5^{-1},4_{II}^{+2},3^{+2}$\\
\hline
     &     &    &        &           &                     &$12\aaa_1$&$19$  &
$8_1^{+1},3^{+2}$\\
\hline
\hline
$65$ &$96$ &$227$&$2^4D_6$& $18$&$2_{II}^{-2},4_7^{+1},8_1^{+1},3^{-1}$&$4\aaa_1$&$19$&
$4_3^{-3},3^{-1}$\\
\hline
     &     &    &        &      &                                      &$8\aaa_1$&$19$&
$2_{II}^{-2},8_7^{+1},3^{-1}$\\
\hline
     &     &    &        &      &                                     & $12\aaa_1$&$19$&
$4_5^{+3}$\\
\hline
     &     &    &        &      &                                     & $16\aaa_1$&$19$&
$2_{II}^{+2},4_3^{-1},3^{-1}$\\
\hline
\hline
 $75$&$192$&$1023$&$4^2\AAA_4$&$18$&$2_{II}^{-2},8_6^{-2}$& $16\aaa_1$&$19$ &$2_{II}^{-2},8_5^{-1}$\\
\hline

\end{tabular}

\end{table}

\vfill


V.V. Nikulin
\par Steklov Mathematical Institute,
\par ul. Gubkina 8, Moscow 117966, GSP-1, Russia;

\vskip5pt

\par Department of Mathematical Sciences,
\par University of Liverpool,
\par Liverpool, L69 3BX, UK

\vskip5pt

\par E-mail: nikulin@mi-ras.ru\ \ \ vnikulin@liv.ac.uk\ \

\end{document}